\documentclass{article}
\usepackage{graphicx}
\usepackage[english]{babel} 
\usepackage{xcolor}
\usepackage{xcolor,colortbl}
\usepackage[document]{ragged2e}
\usepackage{enumitem}
\usepackage{amsmath}
\usepackage{mathtools}
\usepackage{floatrow}
\usepackage{caption}
\usepackage{booktabs}
\usepackage[bottom=1.8cm,top=1.3cm,left=2.6cm,right=2.6cm,includehead,includefoot]{geometry}
\usepackage{abstract}

\usepackage{titlesec} 
\usepackage{fancyhdr} 
\usepackage{titling} 
\usepackage{hyperref} 
\usepackage{amssymb}
\usepackage{xspace}

\usepackage{listings}
\usepackage{listings}
\definecolor{Gray}{gray}{0.95}
\title{\large\bfseries Generalized solutions of Polynomial and Transcendental Equations\\ by the strong method of «Generalized Iterative Method approximation of Roots (GRIM)»\\ \vspace{\baselineskip}}
\author{N.Mantzakouras \& Biragova Eteri\\ 
\normalsize \href{mailto:nikmatza@gmail.com}{e-mail: nikmatza@gmail.com}\\
\normalsize Athens – Greece}
\date{}
\renewcommand{\maketitlehookd}{
\begin{abstract}
\noindent 
\justifying
\color{gray}
In this paper, we explain a new Iterative Method - Fixed Point and develop its convergence theory for finding approximate solutions of nonlinear equations in the setting of Banach spaces. First, we discuss the convergence analysis of our method by separating the equation into functions from which it is derived and the remaining part of the equation, according to the Generalized Theorem[1] for solving polynomial and transcendental equations. Solving an equation, either polynomial or transcendental, is solved only by categorizing the roots and not by some procedure of searching in individually intervals and this way was followed in the past by all existing methods. But the methods developed in the past are limited to isolated intervals without general acceptance. Finally, the method is strengthened by Newton's method to speed up the finding of the roots. At the end of the analysis we give several examples with application of the method.\\
\begin{flushleft}
\textbf{Keywords:} nonlinear equation, polynomial and transcendental equations, iterative Fixed Point Method, continuity conditions, Newton type method, initial values of roots.
\end{flushleft}
\end{abstract}
}
\begin{document}
\maketitle
\justifying
\normalsize
\begin{flushleft}
\end{flushleft}
\textbf{Part I.}\\\\
\textbf{Introduction. "Definitions \& basically Theorems."}\\\\
\textbf{I.1. Inverse Function.} The notion of an inverse is used for many types of mathematical constructions. For example, if we define $f: T \rightarrow S$ is a function restricted to a domain $\mathrm{S}$ and also a range $\mathrm{T}$ in which it is bijective and $g: S \rightarrow T$ is a function satisfying $f(g(s))=s$ for all $s \in S$, then $g$ is the unique function with this property, called \textbf{inverse Function}. It also follows that for all, so, i.e., inversion is symmetric. However, "inverse functions" are also commonly defined for functions that are not bijective (most commonly for elementary functions in the complex plane, which are ), in which case, one of both of the properties may fail to hold (only for genuinely increase function).
$$\text{Apply the forms}
\def\arraystretch{1.5}\left\{\begin{array}{l}
f\left(f^{-1}(y)\right)=y \& y \in f(T) \\
f^{-1}(f(x))=x \& x \in T \\
f^{-1}\left(f^{-1}(x)\right)=f(x)
\end{array}\right\}$$\\
Inverses are also defined for elements of groups, rings, and fields (the latter two of which can possess two different types of inverses known as additive and multiplicative inverses). Every definition of inverse is symmetric and returns the starting value when applied twice.
\newpage\noindent
\textbf{I.2. Fundamental Algebra Theorem.}\\
If $f(x) \in F(x)$ and $\operatorname{deg} \mathrm{f}(\mathrm{x})>0$, then there exists $a \in C$ such that $\mathrm{f}(\mathrm{a})=0$. That is, any non-constant polynomial with coefficients from the set of complex numbers have at least one root in C.\\\\
\textbf{I.3. Definition.I.} If $f(x) \in F(x)$ where $F$ field, with $\operatorname{deg} \mathrm{f}(\mathrm{x})=\mathrm{n}$, is analyzed in product of linear factors (splits) in $\mathrm{F}[\mathrm{x}]$ if $f(x)=c\left(x-a_{1}\right) \cdot\left(x-a_{2}\right) \cdots\left(x-a_{n}\right)$ where $c, a_{1}, a_{2}, \ldots a_{n} \in F$.\\\\
\textbf{I.4. Definition.II.} Let $\mathrm{E} / \mathrm{F}$ be a field extension and $a \in E$. The $\alpha$ is called algebraic above $\mathrm{F}$ if $f(x) \in F(x)$ so that $\mathrm{f}(\mathrm{x})=0$ and $\mathrm{f}(\boldsymbol{\alpha})=0$. If $\boldsymbol{\alpha}$ is not algebraic over $\mathrm{F}$, then $\boldsymbol{\alpha}$ is called transcendental by F.An extension E of $F$ such that the polynomial $f(x)$, when we call it as the polynomial of $E[x]$, has all its roots in $\mathrm{E}$, i.e. analyzed in a product of a linear product to $\mathrm{E}[\mathrm{x}]$. When we talk about fields, $\mathrm{E}$ is a field extension of roots.\\\\
\textbf{I.5. How they are defined as Radicals \& Periodic Radicals as Roots of $\boldsymbol{f(x)=0}$.} If polynomial $f(x) \in C(x)$ with $\operatorname{deg}(\mathrm{x})<5$ has a field of roots $\mathrm{E}$, it has as any root, a sum by finite count of radicals. If $\operatorname{deg} f(x)>4$ except for the cases analyzed in a product of linear factors, the roots take the form of sum of infinite radicals.Also as roots we have Sum of infinite radicals and for transcendental equations as we' ll see above.\\\\
\textbf{I.6. Significant Definition III." Categorization of roots in transcendental equations."}\\
To find roots in transcedental equations we do not follow the procedure as in polynomials down to the fourth degree or Galois theory or other special method. We \textbf{categorize each mononym function} from which consists of the equation and after that we find the corresponding \textbf{inverse function}.\\\\
The corresponding theorem is on page 5 and the formulas for the known functions on page 12. Of course and for \textbf{polynomials up to 4th degree} we use the same theorem.\\\\
\textbf{I.7. Cyclotomic Polynomials.} Consider the cyclotomic polynomial $x^{n}-a=0$. We know that a polynomial is solvable if and only if its Galois group is solvable. There is a general method (we will explain in more detail below) for solving polynomial equations but other cases more complex with this solution of the equation. And this fact is because it involves the highest degree n. In the field of complex numbers, the solution of cyclotomic equation where $\mathrm{n}$ is an integer, are given by DeMoivre's theorem:
\[x_{k}=\sqrt[n]{a} \cdot \zeta_{k}=\sqrt[n]{a} \cdot e^{2 \pi i k / n} \text { where } \mathrm{k}=0,1,2,3 \ldots \mathrm{n}-1.\]
\textbf{I.8. Iteration Method.} In computational mathematics, an \textbf{iterative method }is a mathematical procedure that uses an initial value to generate a sequence of improving approximate solutions for a class of problems, in which the $n$-th approximation is derived from the previous ones. A specific implementation of an iterative method, including the termination criteria, is an algorithm of the iterative method. An iterative method is called \textbf{convergent} if the corresponding sequence converges for given initial approximations. A mathematically rigorous convergence analysis of an iterative method is usually performed; however, heuristic-based iterative methods are also common. Iterative methods are often the only choice for nonlinear equations. However, iterative methods are often useful even for linear problems involving many variables (sometimes of the order of millions), where direct methods would be prohibitively expensive (and in some cases impossible) even with the best available computing power. \\\\
\textbf{I.9. Attractive fixed points.}\\
If an equation can be put into the form $f(x)=x$, and a solution $\mathbf{x}$ is an attractive fixed point of the function $f$, then one may begin with a point $x_{1}$ in the basin of attraction of $\mathbf{x}$, and let $x_{n+1}=f\left(x_{n}\right)$ for $n \geq 1$, and the sequence $\left\{x_{n}\right\}_{n \geq 1}$ will converge to the solution $\mathbf{x}$. Here $x_{n}$ is the $n$th approximation or iteration of $x$ and $x_{n+1}$ is the next or $n+1$ iteration of $x$. Alternately, superscripts in parentheses are often used in numerical methods, so as not to interfere with subscripts with other meanings.\\\\
For example, $x^{(n+1)}=f\left(x^{(n)}\right)$. If the function $f$ is continuously differentiable, a sufficient condition for convergence is that the spectral radius of the derivative is strictly bounded by one in a neighborhood of the fixed point. If this condition holds at the fixed point, then a sufficiently small neighborhood (basin of attraction) must exist.\\\\
\textbf{I.10. Definition IV (Polynomials).} Let $F \subset C$ be a field and $x$ element of $C$, algebraic above it $F$. We will say that element $\mathrm{x}$ is expressed with radicals if it belongs to a field $\mathrm{E}$ such that there is a sequence of fields\\
\[F=F_{0}\subset F_{1}\subset F_{2}\subset \hdots \subset F_{i}\subset F_{i+1}\subset \hdots \subset F_{s}\subset E\]\\
for some integers $s$, where if we assume $x_{0}$ an initial value of $x$ we can take\\
\[F_{i+1}=F_{i}\left(\sqrt[n]{x_{0}}\right)\]\\
for some $x_{0} \in F_{i}$ and $\mathrm{n}$ is Natural number, $0 \leq i \leq s-1$. With $\sqrt[n]{x_{0}}$ we denote a root of polynomial $x^{n}-x_{0}=0 \in F_{i}[x]$. The sequence is called a radical sequence of fields, and the $\mathrm{E} / \mathrm{F}$ extension\textbf{ is called a radical extension}. A polynomial $f(x) \in F[x]$ is called \textbf{resolved by radicals above $\mathbf{F}$, if any it is expressed in radically}, that is, if any a radical extension containing the analysis of $f(x)$ above $F$.\\\\
\textbf{I.11. Galois theorem.} Let $F \subset C$ be a body, $f(x) \in F[x]$ with a field L. The polynomial $\mathrm{f}(\mathrm{x})$ is solvable above $Q$ although only if the Gal (L/F) group is solvable. (At this point we have t\textbf{o give an explanation}. We will not consider here the criteria for an extension L to be solved, so we're not interested in the fact that there is a symmetric polynomial that consider Galois theorem. We will solve by repeating each set of roots in a general polynomial form. So an extension of our own is not necessary to be resolved radically by Galois theory. Therefore, our extension $L^{\prime}$ by radicals may not be solved by Galois theory but \textbf{by the method of repetition} that we accept).\\[3\baselineskip]
\textbf{Part II.}\\\\
\textbf{I.1. Theorem 1.Generalized Theorem of existence and of finding Roots to a Equation.(G$_{\text{T}}$RE)}\\ 
We assume that we have general, elementary functions (in general transcendental) that are analytic except for some isolated singularities and branch cuts, in which case these and their local inversions will have
convergent Taylor series expansions on suitable disks. Under these conditions i consider it a transcendental $f: C \rightarrow C$ where $f(x)=\sum\limits_{i=1}^{n} a_{i} \sigma_{i}(x)+a_{0}, a_{i}, a_{0} \in C$ with $\sigma_{i}(x)$ \textbf{functional factors} where $\sigma_{i}: C \rightarrow C$ for which there is likewise an inversion and if suppose $G$ the field of the roots of equation $f(x)=0$, that defined from relation by category of functional factors $\sigma_{i}$:\\
\[G=F\left(G_{1}=\left\{r_{1}^{\sigma_{1}}, r_{2}^{\sigma_{1}}, \ldots, r_{m_{1}}^{\sigma_{1}}\right\}, G_{2}=\left\{r_{1}^{\sigma_{2}}, r_{2}^{\sigma_{2}}, \ldots, r_{m_{2}}^{\sigma_{2}}\right\},\ldots , G_{k}=\left\{r_{1}^{\sigma_{k}}, r_{2}^{\sigma_{k}}, \ldots, r_{m_{k}}^{\sigma_{k}}\right\}\right) \ldots,\]
\[\left.G_{n}=\left\{r_{1}^{\sigma_{n}}, r_{2}^{\sigma_{n}}, \ldots, r_{m_{n}}^{\sigma_{n}}\right\}\right)\]\\
with $1 \leq k \leq n$ and $0 \leq m_{k} \leq \infty$ where $G_{i}$ \textbf{subfields} of total field of roots where produced the $\sigma_{i}$, and $G$ it will be at the same time and the total global field, which belong the roots of the equation $f(x)=0$.\\\\
If we define the relation of function and divide it into 2 parts ie.\\ \[\sigma_{k}(x)=u_{k}=\frac{1}{a_{k}}\left(-\sum_{i=1, i \neq k}^{n} a_{i} \sigma_{i}(x)-a_{0}\right)\]\\
then the subfields of roots that produces the $\sigma_{k}(x), \ 1 \leq k \leq n$ (according to the generalized theorem by Lagrnage method and the (part\#1(pag.1) that we mentioned [1]), and they will correspond to $G_{k}$. Now if we accept all the functions $\sigma_{k}(x)$ with $1 \leq k \leq n$ and exist inverse functions (invertible procedure) i.e are 1-1, the equation $f(x)=0$ is fully resolved, and is represented by the field $G$. The counting of Roots of equation will be given from the relation $|G|=\left|G_{1} \cup G_{2} \cup \ldots \cup G_{n}\right|$ where $\left(\mathrm{G}_{1}, \mathrm{G}_{2}, \ldots \mathrm{G}_{\mathrm{n}}\right)$ the individual ones subfields) and will be a bounded integer or an infinite integer.\\\\\\
\textbf{Proof.}\\\\
According to the equation $ \sum\limits_{i=1}^{n} a_{i} \sigma_{i}(x)+a_{0}=0, a_{i}, a_{0} \in C, i \in N^{+}$ apply the relations:\\

$$\begin{aligned}
&\hspace{2ex}\sigma_{1}(x)=u_{1}=\frac{1}{a_{1}}\left(-\sum_{i>1}^{n} a_{i} \sigma_{i}(x)-a_{0}\right) \\[10pt]
&\hspace{2ex}\sigma_{2}(x)=u_{2}=\frac{1}{a_{2}}\left(-\sum_{i=1, i \neq 2}^{n} a_{i} \sigma_{i}(x)-a_{0}\right)\\[5pt] 
(I) &\hspace{2ex}.......................................................................\\[5pt]
&\hspace{2ex}\sigma_{k}(x)=u_{k}=\frac{1}{\alpha_{k}} \left(i-\sum_{i=1, i \neq k}^{n} \alpha_{i} \sigma_{i}(x)-\alpha_{0}\right)\\[5pt]
&\hspace{2ex}.......................................................................\\[5pt]
&\hspace{2ex}\sigma_{n}(x)=u_{n}=\frac{1}{\alpha_{n}}\left(-\sum_{i=1, i \neq n}^{n-1} \alpha_{i} \sigma_{i}(x)-\alpha_{0}\right)
\end{aligned}$$\\[4pt]
\noindent
From relations (I) and the property of inversion and with the iterative procedure applying from the previous ones relations we will takes\\
$$\begin{aligned}
&\hspace{2ex}{ }_{\sigma_{1}}^{m_{1}} x_{k_{q}}^{\mu}=\sigma_{1}^{-1}\left(\frac{1}{\alpha_{1}}\left(-\sum_{i>1}^{n} \alpha_{i} \sigma_{i}\left(x_{j}\right)-\alpha_{0}\right)\right)=\sigma_{1}^{-1}\left(u_{1}\right) \in G_{1} \\[12pt]
&\hspace{2ex}{ }_{\sigma_{2}}^{m_{2}} x_{k_{q}}^{\mu}=\sigma_{2}^{-1}\left(\frac{1}{a_{2}}\left(-\sum_{i=1,1 \neq 2}^{n} a_{i} \sigma_{i}\left(x_{t}\right)-a_{0}\right)\right)=\sigma_{2}^{-1}\left(u_{2}\right) \in G_{2}\\[7pt]
(II) &\hspace{2ex}..........................................................................................................\\[7pt]
&\hspace{2ex}{ }_{\sigma_{k}}^{m_{k}} x_{k_{q}}^{\mu}=\sigma_{k}^{-1}\left(\frac{1}{\alpha_{k}}\left(-\sum_{i=1, i \neq k}^{n} \alpha_{i} \sigma_{i}\left(x_{j}\right)-\alpha_{0}\right)\right)=\sigma_{k}^{-1}\left(u_{k}\right) \in G_{k}\\[7pt]
&\hspace{2ex}..........................................................................................................\\[7pt]
&\hspace{2ex}{ }_{\sigma_{n}}^{m_{n}} x_{k_{q}}^{\mu}=\sigma_{n}^{-1}\left(\frac{1}{\alpha_{n}}\left(-\sum_{i=1}^{n-1} \alpha_{i} \sigma_{i}\left(x_{j}\right)-\alpha_{0}\right)\right)=\sigma_{n}^{-1}\left(u_{n}\right) \in G_{n}
\end{aligned}$$\\
\noindent
with $\mu$ the number of repeats $, 1 \leq k \leq n, \ 0 \leq m_{k} \leq \infty, \  1 \leq \mathrm{q} \leq m_{k}$ which means that the count of roots it can to be found in interval $[0, \infty)$ and apply in generally for number $\gamma_{m_{k}}$ of the coefficient $m_{k}$ (for any subfield $G_{k}$ number of roots) will be $\gamma_{m_{k}} \in \mathrm{N}^{+}$. The transformations of relation (II) they produce exactly the subfields $G_{1}, G_{2}, \ldots, G_{k}, \ldots G_{n}$ which generally they constitute the $\mathrm{G}$. The genitive form finding of roots of $f(x)=\sum\limits_{i=1}^{n} a_{i} \sigma_{i}(x)+a_{0}=0, \ \{x,a_{i}, a_{0}\in \mathbb{C},i\in\mathbb{N}^{+} \}$ with the Periodic Radicals will be done with the generalized function $\Phi_{k}\left(u_{k}\right)=\sigma_{k}^{-1}\left(u_{k}\right)$ resulting from every relations in simple $\sigma_{k}(x)=u_{k} \Rightarrow x_{k_{q}}^{\sigma_{k}}=\sigma_{k}^{-1}\left(u_{k}\right)$ where $\sigma_{k}(x)$ (\textbf{functional factors of equation} $1\leq k\leq n$) and each root will given from \\

$$_{\sigma_{k}}^{m_{k}}x_{k_{q}}^{\mu}=\Phi_{k,u}\left(\frac{1}{\alpha_{k}}\left(-\sum\limits_{i=1, i\neq k}^{n}\alpha_{i}\sigma_{i}(\#)-\alpha_{0}\right)\right) $$ \\
\noindent\\
where the symbol \# shows repetitive (periodic) procedure of values. For \textbf{a number of repetitions} we find each time the roots $m_{k}$, $\mu$ the number of repeats, of category $\mathrm{k}$ from the form $_{\sigma_{k}}^{m_{k}}x_{j}^{\mu}$ until some final value, with approach of error $<10^{\wedge}(-\rho)$, $\rho$ is integer and $\rho>>1$. The following stopping criteria is used for computer programs: 1.) $\left|x_{j+1}-x_{j}\right|<\varepsilon$, 2.) $\left|f\left(x_{j}\right)\right|<\varepsilon$ where $\varepsilon=10^{\wedge}(-\rho), \  \rho>1$, all Integers. Repeat procedures for all subfields of roots, are expressed computationally by the relations that result from (II) forms with programming language in the form (III). Analytically we will have them relations\\
\[\begin{aligned}
&{ }_{\sigma_{1}}^{m_{1}} x_{1_{q}}^{\mu}=\Phi_{1, u_{1}}\left(\frac{1}{\alpha_{k}}\left(-\sum_{i>1}^{n} \alpha_{i} \sigma_{i}(\#)-\alpha_{0}\right)\right),\\\\
&{ }_{\sigma_{2}}^{m_{2}} x_{2_q}^{\mu}=\Phi_{2, u_{2}}\left(\frac{1}{a_{2}}\left(-\sum_{i=1, i \neq 2}^{n} a_{i} \sigma_{i}(\#)-\alpha_{0}\right)\right), 
\end{aligned}\]\\\\
${ }_{\sigma_{n}}^{m_{n}} x_{n_{q}}^{\mu}=\Phi_{n, u_{n}}\left(\frac{1}{\alpha_{n}}\left(-\displaystyle \sum_{i=1, i\neq n}^{n} \alpha_{i} \sigma_{i}(\#)-\alpha_{0}\right)\right)$ with $\gamma_{m_{k}} \in N^{+}, k \in\{1,2 \ldots n\}, \mu$ the number of repeats.\\\\
Obviously if we have a limited number of functions which to be $1-1$ and be inverses functions we will limit ourselves to the subfields they produce as many they can. The counting of Roots of equation will be given from the relation $|G|=\left|G_{1} \cup G_{2} \cup \ldots \cup G_{n}\right|$ where $G_{1}, G_{2}, \ldots G_{n}$ (the individual ones subfields) where for every one $\mathrm{G}_{\mathrm{k}}$ corresponds an one inverse $\sigma_{k}^{-1}$ and the final count of roots of solutions of ${f}(x)=0$ results from the sum:\\
\[|\mathrm{G}|=\sum_{\mathrm{i}=1}^{\mathrm{n}} \sum_{\mathrm{j}=1}^{m_{j}} \mathrm{q}_{\mathrm{j}}^{\mathrm{i}} \text{ (Total number - solutions)}\]
also $q_{j}^{i}$ the number of roots per $m_{i}$ category in the i position respectively (any category $m_{i}$ it has subcategories with count number from 1 to $\mathrm{m}_{\mathrm{i}}$) and finally the $\gamma_{\mathrm{m}_{\mathrm{q}}}$ is the number of subcategories (per i) category.\\\\\\
\textbf{II.2. Subfields, categories of subfields and number of roots per category [1]}\\\\
The $f(z)=\sum\limits_{i=1}^{n} \alpha_{i} \cdot \sigma_{i}(z)-\alpha_{0}=0,\{z, \alpha_i, \alpha_{0}\in\mathbb{C},i\in\mathbb{N}^{+}\}$ it has roots $Z_{s_{q}}^{L_{k, q}}$, where $s_{q} \in Z, q \in N$ and $s_{q}$ is a multiple parameter with $q$ specifying the number of categories within the corpus itself associated with either complex exponential or trigonometric functions. The $L_{k, q}$ \textbf{is the subfield} concerns only this form, that is to say $\left(\sigma_{k}\right)$. But the \textbf{k-subfield itself} can have \textbf{$\mathbf{q}$ categories}. In case we want to connect the subfields of the root in the complex plane and the general solution equation per category we will have the general relation:\\\\
\[_{\sigma_{k}}^{L_{k,q}}z_{s_{q}}^{\mu}=\sigma_{k}^{-1}\left(\frac{\alpha_{0}}{\alpha_{k}}-\frac{1}{\alpha_{k}} \sigma_{k}^{c}\left(\sigma_{k}^{-1}\left(\frac{\alpha_{0}}{\alpha_{k}}-\frac{1}{\alpha_{k}} \sigma_{k}^{c}\left(\sigma_{k}^{-1}\left(\frac{\alpha_{0}}{\alpha_{k}}-\ldots-\right.\right.\right.\right.\right.\]
\[\left.\left.\left.\left.\left.-\frac{1}{\alpha_{k}} \sigma_{k}^{c}\left(\sigma_{k}^{-1}\left(\frac{\alpha_{0}}{\alpha_{k}}-\frac{1}{\alpha_{k}} \cdot \sum_{\mathrm{i}=2, \mathrm{i} \neq \mathrm{k}}^{\mathrm{n}} \alpha_{\mathrm{i}} \sigma_{\mathrm{i}}\left(\mathrm{x}_{0}\right)\right)\right)\right)\right)\right)\right)\right)\tag{1$^{\ast}$}\]\\
\[1 \leq k \leq n, 0 \leq s \leq \infty, q \in N^{+} \text{ where } \sigma_r^c \text{ the complement of } \sigma_r.\]\\
\noindent
Now, for the generalisation of cases, because this $k$ takes values from to $1 \div n$, consequently \textbf{the count of the basic subfields} of roots also will be $n$, and consequently \textbf{the field of total of the roots} of the equation \textbf{is L} and will apply:\\
\[
L=\bigcup_{k=1}^{n} L_{k}^{u_{k}}, \quad L_{k}^{u_{k}}=\bigcup_{q=1}^{u_{k}} L_{k, q},\left\{k, u_{q} \in N, 1 \leq k \leq n\right\}\tag{2$^{\ast}$}
\]\\
\noindent
The parameter $u_{q}$ takes different values and depends on whether a term function of the equation is trigonometric or exponential or multiple exponential or polynomial. Example for trigonometric is $1 \leq u_{q} \leq 2$ and exponential with multiplicity 1 degree is $u_{q}=1$. We will look at these specifically in examples below. $s_{q}$ is a parameter with $q$ specifying the number of categories within the body itself that are associated with either exponential or trigonometric functions or polynomial. "But what is enormous interest about formula (1$^{\ast}$) is that it gives us all the roots for each term of the equation $\sigma(z)=0$, once it has been analyzed and categorized using the inverse function technique $\sigma_{k}^{-1}(x_{0},s_{q})$ after we do the analysis separately in each case.\\\\

\noindent\textbf{II.3. Number of Roots by use sets.} In its general form, the principle of inclusion-exclusion states [2] that for finite sets $A_{1}, \ldots, A_{n}$, one has the identity and just let $A_{i}$ be the subset of elements of $S$. Then we have
\[\begin{aligned}
&\left|A_{1} \cup A_{2} \cup \cdots \cup A_{n}\right|= \sum_{1 \leq i \leq n}\left|A_{i}\right|-\sum_{1 \leq i<j \leq n}\left|A_{i} \cap A_{j}\right| \\
&+\sum_{1 \leq i<j<k \leq n}\left|A_{i} \cap A_{j} \cap A_{k}\right|-\cdots+(-1)^{n+1}\left|A_{1} \cap A_{2} \cap \cdots A_{n}\right| \\
&\left|A_{1} \cup A_{2} \cup \cdots \cup A_{n}\right|=\sum_{k=1}^{n}(-1)^{k+1} \sum_{J \subseteq[n]:|J|=k}\left|\cap_{j \in J} A_{j}\right|
\end{aligned}\]\\
Here we have the equality $G_{i}=A_{i} \forall i \in\{1,2 \ldots n\}$ and $G=S$. Clearly the subfields can be infinite number elements. The $G_{i}$ are subsets of the $G$ set of Roots of the equation $f(x)=0$.\\
\textbf{Number intersections.} The proof is done by induction:\\
\[C\left(\sum_{1 \leq i \leq n}\left|A_{i}\right|\right)=c_{n}^{1}=n,\text{ } C\left(\sum_{1 \leq i \leq n}\left|A_{i} \cap A_{J}\right|\right)=c_{n}^{2}=n \cdot(n+1) / 2, \ldots, C\left(\left|A_{i} \cap A_{J} \cap \ldots\right|\right)=c_{n}^{n}=1\]\\\\
\textbf{The total count} $=2^{n}-1$.
\newpage\noindent
\textbf{Example:} The principle is more clearly seen in the case of three sets $(n=3)$, which for the sets $A, B$ and $C$ is given by
\[S=|A \cup B \cup C|=|A|+|B|+|C|-|A \cap B|-|A \cap C|-|B \cap C|+|A \cap B \cap C|\]\\
\begin{figure}[h!]
  \centering
    \includegraphics[width=0.3\textwidth]{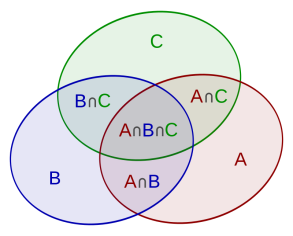}
\caption{} \label{f4.defins0}
\end{figure}\\\noindent
This formula can be verified by counting how many times each region in the Venn diagram (fig.1) figure is included in the right-hand side of the formula. In this case, when removing the contributions of overcounted elements, the number of elements in the mutual intersection of the three sets has been subtracted too often, so must be added back in to get the correct total.
\textbf{This theory applies in our case both to the union of the principal fields of the basic functions and to
the subfields of the partial functions.}\\\\\\
\textbf{Part III. Method Infinity Periodicity of Approximate Roots (GRIM)}\\\\
\textbf{III.1. Introduction.}\\\\
This method applies to polynomials with real or complex coefficients as well as transcedental polynomial formats. For such polynomials, we know that all complex roots appear as coupled pairs. The method attempts to locate the zeros of such polynomials by searching for zeros in a repetitive process after being split into two segments, and then in one section we calculate the inverse function that is the polynomial from mononym to until a trinomial. Therefore, we make the following reasoning, based on the theorem of continuity and if a function it has derivative and inverse.\\\\
\textbf{From special to general method:} Suppose we have the function $f(x)=\displaystyle\sum_{i=0}^{n} \alpha_{i} x^{i}=0$, (1) we divide the equation with a sum of two functions equal to zero. So we will have a monomial and a polynomial. This namely $f_{1}(x)-f_{2}(x)=0$. Because we are looking to find $\mathrm{x}$, it will be the same root for both functions that as a sum are an equation we can equalize the one function $f_{1}(x)=u$ and $f_{2}(x)=u$ then with the relation $x_{k_{1}}=f_{1}^{-1}(u)=f_{1}^{-1}\left(f_{2}\left(x_{m}\right)\right)$ and same time apply $x_{k_{2}}=f_{2}^{-1}(u)=f_{2}^{-1}\left(f_{1}\left(x_{m}\right)\right)$ for $k_2$ count of Roots for the category of function $f_2$. Finally, with periodic iterative radicals we will have the final relation at the root of the initial equation (1) where $x_{k_{2}}^{r}=f_{2}^{-1}\left(f_{1}\left(x_{\#}\right)\right), r=0,1,2, \ldots N$ for $r$ repetitions of maximum number $N$ and certainly in this case we call relation $f_{1}(x)=\displaystyle\sum_{i=0}^{n-1} a_{i} x^{i}, f_{2}(x)=-a_{n} x^{n}$.\\\\\\
According to the Theorem and method \textbf{G.R.L-N[1] the total of the roots will be the Union of the roots} that will be found by the solution for the method \textbf{(GRIM)} for each monomial of polynomial $\displaystyle\sum_{i=0}^{n} a_{i} x^{i}=0$. Therefore, because $i$ have $n$ monomials will be $n$ Roots. It is of course equivalent to the solution of the monomial with the maximum exponent, that is $n$. Of course this does not apply to transcendental equations in general, because there we have another categorization with complex parameters.\\\\\\
\textbf{III.2.Theorem 2. Generalized Iterative Method approximation of Roots (GRIM)}\\\\
We assume that we have general, elementary functions(in general transcendental) that are analytic except for some isolated singularities and branch cuts, in which case these and their local inversions will have convergent Taylor series expansions on suitable disks. Under these conditions for the determination of the roots any equation $f(x)=\sum\limits_{i=1}^{n} a_{i} \sigma_{i}(x)-a_{0}=0, \ \{x,a_{i},a_{0}\in\mathbb{C},i\in\mathbb{N}^{+} \} \ (1)$ where $f: C \rightarrow C$ which consists from distinct monomials functions $\sigma_{i}(x), 1 \leq i \leq n \wedge a_{i} \in C$ and for which there is likewise an inversion, after we accept $\sigma_{0}(x)=1$, exist a repeating relation for every subfield of roots $G_{l}, G_{2}, \ldots, G_{r}, \ldots G_{n}$ and $\mathrm{k}$ the counting of the roots for any $G_{r}, r \in\{1,2, \ldots, n\}$. Now if we now assume that that $x_{0}$ initial value for $x$ then and from GR$_{\text{T}}$E it will hold more closely with indices and for the any roots of $\mathrm{f}(\mathrm{x})$ will apply:\\\\
\[\begin{aligned}
i) \text{ }{ }_{\sigma_{r}}^{m_{r}} x_{k_{r}}^{\mu}&=\sigma_{r}^{-1}\left(\frac{\alpha_{0}}{\alpha_{r}}-\frac{1}{\alpha_{r}} \sigma_{r}^{c}\left(\sigma_{r}^{-1}\left(\frac{\alpha_{0}}{\alpha_{r}}-\frac{1}{\alpha_{r}} \sigma_{r}^{c}\left(\sigma_{r}^{-1}\left(\frac{\alpha_{0}}{\alpha_{r}}-\ldots-\right.\right.\right.\right.\right.\\
&\left.\left.\left.\left.\left.-\frac{1}{\alpha_{r}} \sigma_{r}^{c}\left(\sigma_{r}^{-1}\left(\frac{\alpha_{0}}{\alpha_{r}}-\frac{1}{\alpha_{r}} \cdot \sum_{\mathrm{i}=1, \mathrm{i} \neq \mathrm{r}}^{\mathrm{n}} \alpha_{i} \sigma_{\mathrm{i}}\left(\mathrm{x}_{0}\right)\right)\right)\right)\right)\right)\right)\right), 1<\mathrm{r} \leq \mathrm{n}, 1 \leq \mathrm{k} \leq \infty
\end{aligned}\]\\
\[\begin{aligned}
ii)\text{ } { }_{\sigma_{r}}^{m_{r}} x_{k_{r}}^{\mu}&=\sigma_{1}^{-1}\left(\frac{\alpha_{0}}{\alpha_{1}}-\frac{1}{\alpha_{1}} \sigma_{1}^{c}\left(\sigma_{1}^{-1}\left(\frac{\alpha_{0}}{\alpha_{1}}-\frac{1}{\alpha_{1}} \sigma_{1}^{c}\left(\sigma_{1}^{-1}\left(\frac{\alpha_{0}}{\alpha_{1}}-\ldots-\frac{1}{\alpha_{1}} \sigma_{1}^{c}\left(\sigma_{1}^{-1}\left(\frac{\alpha_{0}}{\alpha_{1}}-\right.\right.\right.\right.\right.\right.\right.\\
&\left.\left.\left.\left.\left.\left.\left.-\frac{1}{\alpha_{1}} \cdot \sum_{\mathrm{i}=2}^{\mathrm{n}} \alpha_{i} \sigma_{\mathrm{i}}\left(\mathrm{x}_{0}\right)\right)\right)\right)\right)\right)\right)\right), \mathrm{r}=1,1 \leq \mathrm{k} \leq \infty
\end{aligned}\]\\
where $\sigma_{r}^{c}$  the complement of $\sigma_{r}$.\\\\
where $\mu$ the number of repetitions and is based on the iterative procedure of generalized theorem GR$_{\text{T}}$E. This relation is also the generalized recurrent periodically algorithm of finding of the roots of any equation.\\\\
\textbf{Proof}\\\\
For proof we go from the special to the general case. Suppose we have the function $f[x]=\displaystyle\sum_{i=1}^{\infty} \alpha_i \sigma_i(x)-\alpha_0=0$ (1) where $f: C \rightarrow C, \alpha_{i} \in C$ where $\alpha_{i}$ constants. I separate the initial equation with a sum of functions and a monomial function and then therefore the two parts they will have sum equal to zero. This means as example: $f_{1}(x)+f_{2}(x)-\alpha_0=0$ where: $f_{1}(x)=a_{1}\sigma_{1}(x) \wedge f_{2}(x)=\displaystyle\sum_{i=2}^{n} a_{i} \sigma_{i}(x) \Rightarrow { }_{\sigma_{1}} x_{j+1}=f_{1}^{-1}\left(\dfrac{a_{0}}{a_{1}}-\dfrac{1}{a_{1}}f_{2}\left(x_{j}\right)\right) \ (3)$ which results from the inverse of the first function with the simultaneous repetitive procedure. The result gives one of roots of subfield $G_{1}$. Because we have a iterative process, it will come up with iterations for the given root type:\\\\
\[\begin{aligned}
{ }_{\sigma_{1}}^{m_{1}} x_{k_{1}}^{\mu}&=\mathrm{f}_{1}^{-1}\left(\frac{\alpha_{0}}{\alpha_{1}}-\frac{1}{\alpha_{1}} \mathrm{f}_{2}\left(\mathrm{f}_{1}^{-1}\left(\frac{\alpha_{0}}{\alpha_{1}}-\frac{1}{\alpha_{1}} \mathrm{f}_{2}\left(\mathrm{f}_{1}^{-1}\left(\frac{\alpha_{0}}{\alpha_{1}}-\ldots-\frac{1}{\alpha_{1}} \mathrm{f}_{2}\left(\mathrm{f}_{1}^{-1}\left(\frac{\alpha_{0}}{\alpha_{1}}-\frac{1}{\alpha_{1}}  \sum_{\mathrm{i}=2}^{\mathrm{n}} \alpha_{\mathrm{i}} \sigma_{\mathrm{i}}\left(\mathrm{x}_{0}\right)\right)\right)\right)\right)\right)\right)\right)
\end{aligned}\]
\newpage\noindent
In general we have for $m_{r}$ solutions corresponding to an equation of the of $\sigma_{r}$  function i.e:\\
\[\begin{aligned}
{ }_{\sigma_{r}}^{m_{r}} x_{k_{r}}^{\mu}
=\mathrm{f}_{\mathrm{r}}^{-1}\left(\frac{\alpha_{0}}{\alpha_{r}}-\frac{1}{\alpha_{r}} \mathrm{f}_{\mathrm{r}}^{\mathrm{c}}\left(\mathrm{f}_{\mathrm{r}}^{-1}\left(\frac{\alpha_{0}}{\alpha_{\mathrm{r}}}-\frac{1}{\alpha_{r}} \mathrm{f}_{\mathrm{r}}^{\mathrm{c}}\left(\mathrm{f}_{\mathrm{r}}^{-1}\left(\frac{\alpha_{0}}{\alpha_{r}}-\ldots-\right.\right.\right.\right.\right.\\
\left.\left.\left.\left.\left.-\frac{1}{\alpha_{r}} \mathrm{f}_{\mathrm{r}}^{\mathrm{c}}\left(\mathrm{f}_{\mathrm{r}}^{-1}\left(\frac{\alpha_{0}}{\alpha_{\mathrm{r}}}-\frac{1}{\alpha_{r}} \sum_{\mathrm{i}=1, \mathrm{i} \neq \mathrm{r}}^{\mathrm{n}} \alpha_{\mathrm{i}} \sigma_{\mathrm{i}}\left(\mathrm{x}_{0}\right)\right)\right)\right)\right)\right)\right)\right. (4)
\end{aligned}\]
where $f_{r}^{c}$ the complement of $f_{r}$.\\\\
Here $x_{0}$ is initial value from $\mathrm{x}$ and which results after from $\mu$ repetitions, take into account the same and the other fields $G_{1}, G_{2}, \ldots, G_{r}, \ldots G_{n}$. If I call $\sigma_{\mathrm{r}}^{\mathrm{c}}(\mathrm{x})= \sum\limits_{\mathrm{i} = 1, \mathrm{i}\neq \mathrm{r}}^{\mathrm{n}} \alpha_{\mathrm{i}} \cdot \sigma_{\mathrm{i}}(\mathrm{x})$ is complement of $\sigma_{\mathrm{r}}(\mathrm{x})$ where $x_{0}$ initial value for $x$ then and from GR$_{\text{T}}$E (Theorem 1) will be apply for the any roots of $f(x)$:\\\\
from (1)\\
\[\alpha_{1} \sigma_{1}(x)=-\left(\sum_{i=2}^{n} \alpha_{i} \sigma_{i}(x)-\alpha_{0}\right) \Leftrightarrow \sigma_{1}(x)=-\left(\frac{1}{\alpha_{1}} \sum_{i=2}^{n} \alpha_{i} \sigma_{i}(x)-\alpha_{0}\right)\]
\[\Leftrightarrow x=\sigma_{1}^{-1}\left(-\left(\frac{1}{\alpha_{1}} \sum_{i=2}^{n} \alpha_{i} \sigma_{i}(x)-\alpha_{0}\right)\right) \tag{5}\]\\
and then apply for repeat with simple procedure $_{\sigma_{r}} x_{j+1}=\sigma_{r}^{-1}\left(-\dfrac{1}{a_{r}} \left( \displaystyle\sum_{i \neq r, i=1}^{n} \alpha_{i} \cdot \sigma_{i}\left(x_{j}\right)- \alpha_0\right) \right)$\hspace*{5ex}(6)\\\\
We can now easily deduce the generalized relationship
\[\begin{aligned}
i) \text{ }{ }_{\sigma_{k}}^{L_{k,q}} x_{s_{q}}^{\mu}&=\sigma_{k}^{-1}\left(\frac{\alpha_{0}}{\alpha_{k}}-\frac{1}{\alpha_{k}} \sigma_{k}^{c}\left(\sigma _ { k } ^ { - 1 } \left(\frac{\alpha_{0}}{\alpha_{k}}-\frac{1}{\alpha_{k}} \sigma_{k}^{c}\left(\sigma _ { k } ^ { - 1 } \left(\frac{\alpha_{0}}{\alpha_{k}}-\ldots-\right.\right.\right.\right.\right.\\
&-\left.\left.\left.\left.\left.\frac{1}{\alpha_{k}} \sigma_{k}^{c}\left(\sigma_{k}^{-1}\left(\frac{\alpha_{0}}{\alpha_{k}}-\frac{1}{\alpha_{k}} \cdot \sum_{i=1, i \neq k}^{n} \alpha_{i} \sigma_{i}\left(x_{0}\right)\right)\right)\right)\right)\right)\right)\right)\\
&\text { if } 1<\mathrm{k} \leq \mathrm{n}, 0 \leq \mathrm{s} \leq \infty, \mathrm{q} \in \mathrm{N}^{+} \text {, }
\end{aligned}\]\\
\[\begin{aligned}
ii)\text{ }{ }_{\sigma_{k}}^{L_{k,q}} x_{s_{q}}^{\mu}&=\sigma_{1}^{-1}\left(\frac{\alpha_{0}}{\alpha_{1}}-\frac{1}{\alpha_{1}} \sigma_{1}^{c}\left(\sigma_{1}^{-1}\left(\frac{\alpha_{0}}{\alpha_{1}}-\frac{1}{\alpha_{1}} \sigma_{1}^{c}\left(\sigma_{1}^{-1}\left(\frac{\alpha_{0}}{\alpha_{1}}-\ldots-\right.\right.\right.\right.\right.\\
-&\left.\left.\left.\left.\left.\frac{1}{\alpha_{1}} \sigma_{1}^{c}\left(\sigma_{1}^{-1}\left(\frac{\alpha_{0}}{\alpha_{1}}-\frac{1}{\alpha_{1}} \cdot \sum_{i=2}^{n} \alpha_{i} \sigma_{i}\left(x_{0}\right)\right)\right)\right)\right)\right)\right)\right)\\
&\text { if } \mathrm{k}=1,0 \leq \mathrm{s} \leq \infty, \mathrm{q} \in \mathrm{N}^{+}\\
\end{aligned}\]
\[\text{Apply}: \alpha_{k} \sigma_{k}(x)+\sigma_{k}^{c}(x)-\alpha_{0}=0 \wedge \sigma_{k}^{c}(x)=\sum_{i=1, i \neq k}^{n} \alpha_{i} \sigma_{i}(x) \text{ where } \sigma_{\mathrm{k}}^{\mathrm{c}} \text{ the complement of }\sigma_{k}.\]
where $\sigma_{r}^{c}$  the complement of $\sigma_{r}$.\\\\
where $\mu$ the number of repetitions, $m_{r}$ emerges from the number produced by categorizing the roots for each $G_{r}, r \in\{1,2, \ldots, n\}$, field and $\mathrm{x}_{0}$ initial value for $\mathrm{x}$ that is based on the iterative procedure of generalized theorem GR$_{\text{T}}$E.\\\\\\\\
\textbf{Part IV. Solved examples with the method (GRIM) [2,3,4,5]}\\\\
\textbf{IV.1.1. Solving of general transcendental trinomial}\\
\[\alpha_{1} z^{r_{1}} \cdot a^{z^{q_{1}}}+\alpha_{2} \cdot z^{r_{2}} \cdot b^{z^{q_{2}}}-\alpha_{0}=0 \ \left(\alpha_{i}, r_{i}, q_{i}, a, b \in C, i=1,2 \wedge \alpha_{0} \in C\right)\]\\
\textbf{IV.1.2.} We have the general equation with relation $f(z)=\sum\limits_{i=1}^{2} \alpha_{i} \cdot \sigma_{i}(z)-\alpha_{0}=0,\left\{r_{i}, q_{i}, \alpha_{i}, \alpha_{0} \in C, i=1,2\right\}$. Of course in this case the following will apply: $\sigma_{1}(z)=z^{r_{1}} \cdot a^{z^{q_{1}}} \wedge \sigma_{2}(z)=z^{r_{2}} \cdot b^{z^{q_{2}}}$. But because we have identical monomials we can simply consider the one case it is enough to choose one of the 2 functions, in this case the 1 nd function i.e $\sigma_{1}(z)=z^{r_{1}} \cdot a^{z^{q_{i}}}, r_{1}, q_{1} \in C$. According to the theory, we can to map the \textbf{first term function} to a variable, suppose $u$ and will apply: \[\sigma_{1}(z)=z^{r_{1}} \cdot a^{z^{q_{1}}}=u \Leftrightarrow r_{1} \log z+z^{q_{1}} \log a=\log u+2 k \pi I, q_{1} \neq 0, k \in Z\tag{1}\]
From $(1) \Leftrightarrow\left(\frac{r_{1}}{q_{1}}\right) \cdot \log z^{q_{1}}+z^{q_{1}} \log a=\log u+2 k \pi I(2),\ I:$ the imaginary unit.\\\\
If I call $z^{q_{1}}=y$ (3) then (2) becomes
\[\left(\frac{r_{1}}{q_{1}}\right) \cdot \log y+y \log a=\log u+2 k \pi I.\tag{4}\]
The point here is to relate variable $y$ to a known solvable function as equation. This function is the W-function. If we now call the relationship
\[t y \cdot e^{t y}=s \ (5) \Leftrightarrow t y=W\left(k_{1}, s\right) \Leftrightarrow y=\frac{W\left(k_{1}, s\right)}{t} \ (6), \ \left\{s, t \in C, t \neq 0, k_{1} \in Z\right\}\]
If I now logarithmize relationship (5), I obtain the relationship
\[t y \cdot e^{t y}=s \Leftrightarrow \log y+t y=\log \frac{s}{t}+2 k \pi I\tag{7}\]
Combining $(4,7)$ we arrive at the relationship
\[
t=\frac{\log a}{\frac{r_{1}}{q_{1}}} \ (8) \wedge u^{\frac{q_{1}}{r_{1}}}=\frac{s}{t} \ (9)
\]
By combining relations $(3,6,8,9)$ we obtain the final solution which will be of the form
\[
z_{1}=\left\{\left(\frac{r_{1}}{q_{1} \cdot \log a}\right) \cdot W\left(k_{1}, \frac{\log a}{\frac{r_{1}}{q_{11}}} \cdot u^{\frac{q_{1}}{r_{1}}}\right)\right\}^{\frac{1}{q_{1}}} \cdot e^{\frac{2 k_{2} \pi \cdot \mathrm{I}}{q_{1}}}(10),\left\{k_{2}, k_{1} \in Z\right\}, k_{2}=1 \div q_{1}-1
\]
In the same way we calculate for the second kernel the relation
\[
z_{2}=\left\{\left(\frac{r_{2}}{q_{2} \cdot \log b}\right) \cdot W\left(k_{1}^{\prime}, \frac{\log b}{\frac{r_{2}}{q_{2}}} \cdot u^{\frac{q_{2}}{r 2}}\right)\right\}^{\frac{1}{q_{2}}} \cdot e^{\frac{2 k_{2}^{\prime} \pi \cdot \mathrm{I}}{q_{2}}} \ \ (10),\left\{k_{1}^{\prime}, k_{2}^{\prime} \in Z\right\}, k_{2}^{\prime}=1 \div q_{2}-1
\]
Up to this point we have developed the concept of the iteration process. But the basic GR$_{\text{T}}$E theorem shows that to complete the solution we will have 2 formations for each case separately. i.e $\sigma_{1}(z)=z^{r_{1}} \cdot a^{z^{q_{1}}} \wedge \sigma_{2}(z)=z^{r_{2}} \cdot b^{z^{q_{2}}}$.\\\\
According to the theory GR$_{\text{T}}$E we developed before we consider 2 transformations. Here we have two functions $\sigma_{1}(z), \sigma_{2}(z)$ i.e $\sigma_{1}(z)=z^{r_{1}} \cdot a^{z^{q_{1}}} \wedge \sigma_{2}(z)=z^{r_{2}} \cdot b^{z^{q_{2}}}$. For each case we need to find the inverse function separately per function. Therefore we expect to have 2 subfields of roots $L_{1}, L_{2}$ and therefore the total solution of the equation will be $L=L_{1} \cup L_{2}$. [1]\\\\\\
\textbf{IV.1.3. Finding the $\boldsymbol{L_{1}}$ subfield.}\\\\
The first roots of sub-field $L_{1}$ results from inverse function and we give for kernel the relation
\[\sigma_{1}(z)=u \Rightarrow
z_{\left\{s_{1}, s_{2}\right\}}^{L_{1}}=\left\{\left(\frac{r_{1}}{q_{1} \cdot \log a}\right) \cdot W\left(s_{1}, \frac{\log a}{\frac{r_{1}}{q_{1}}} \cdot u^{\frac{q_{1}}{r_{1}}}\right)\right\}^{\frac{1}{q_{1}}} \cdot e^{\frac{2 s_{2} \pi \cdot \mathrm{I}}{q_{1}}},\left\{s_{2}, s_{1} \in Z\right\}, s_{2}=0 \div q_{1}-1\tag{10}
\]
So, using this method (GRIP) we will have 1 relation to find the set of solutions of the first subfield. We talk about $\boldsymbol{q_{1}}$ \textbf{Categories of infinite (by parameter $\boldsymbol{s_{1}}$)} complex roots.\\\\\\
\textbf{IV.1.4. Finding the $\boldsymbol{L_{2}}$ subfield.}\\\\
The second roots of sub-field $L_{2}$ results from inverse function and we give for kernel the relation 
\[\sigma_{2}(z)=u \Rightarrow
z_{\left\{s_{1}^{\prime}, s_{2}^{\prime}\right\}}^{L_{2}}=\left\{\left(\frac{r_{2}}{q_{2} \cdot \log b}\right) \cdot W\left(s_{1}^{\prime}, \frac{\log b}{\frac{r_{2}}{q_{2}}} \cdot u^{\frac{q_{2}}{r_{2}}}\right)\right\}^{\frac{1}{q_{2}}} \cdot e^{\frac{2 s_{2}^{\prime} \pi \cdot \mathrm{I}}{q_{2}}},\left\{s_{2}^{\prime}, s_{1}^{\prime} \in Z\right\}, s_{2}^{\prime}=0 \div q_{2}-1\tag{11}\]
So, using this method (GRIP) we will have 1 relation to find the set of solutions of the second subfield. We talk about $\boldsymbol{q_{2}}$ \textbf{Categories of infinite (by parameter $\boldsymbol{s_{1}^{\prime}}$)} complex roots.\\\\\\
\textbf{IV.1.5. Iterative Method approximation of Roots.}\\\\
In our case we have 2 functions, with one Sub-fields per function, which means that for the set of solutions of the field $L$ of the equation we have $L=L_{1} \cup L_{2}$. The iterative process that applies according to the theory we have developed will therefore be summarised in the formula\\\\
i) Roots set for for the sub-field $L_{1}$ that generated by the function $\sigma_{1}(z)$\\
\[{ }_{\sigma_{1}}^{L_{1,q}} z_{s_{q}}^{\mu}=\sigma_{1}^{-1}\left(\frac{\alpha_{0}}{\alpha_{1}}-\frac{1}{\alpha_{1}} \sigma_{1}^{c}\left(\sigma_{1}^{-1}\left(\frac{\alpha_{0}}{\alpha_{1}}-\frac{1}{\alpha_{1}} \sigma_{1}^{c}\left(\sigma_{1}^{-1}\left(\frac{\alpha_{0}}{\alpha_{1}}-\ldots- \right.\right.\right.\right.\right.\]
\[\left.\left.\left.\left.\left.-\frac{1}{\alpha_{1}} \sigma_{1}^{c}\left(\sigma_{1}^{-1}\left(\frac{\alpha_{0}}{\alpha_{1}}-\frac{1}{\alpha_{1}} \cdot \sum\limits_{\mathrm{i}=2}^{2} \alpha_{\mathrm{i}} \sigma_{\mathrm{i}}\left(\mathrm{z}_{0}\right)\right)\right)\right)\right)\right)\right)\right), 1 \leq k \leq 2, \ 0\leq s\leq \infty , q\in \mathbb{N}^{+}, \]\\\\
Apply: $\alpha_{1} \sigma_{1}(z)+\sigma_{1}^{c}(z)-\alpha_{0}=0 \wedge \sigma_{1}^{c}(z)=\sum\limits_{i=2}^{2} \alpha_{2} \sigma_{i}(z)=\alpha_{2} \cdot z^{r_{2}} \cdot b^{z^{q_{2}}}$
\[
\sigma_{1}^{-1}(z)=\left\{\left(\frac{r_{1}}{q_{1} \cdot \log a}\right) \cdot W\left(s_{1}, \frac{\log a}{\frac{r_{1}}{q_{1}}} \cdot u^{\frac{q_{1}}{r_{1}}}\right)\right\}^{\frac{1}{q_{1}}} \cdot e^{\frac{2 s_{2} \pi \cdot I}{q_{1}}},\left\{s_{2}, s_{1} \in Z\right\}, s_{2}=0 \div q_{1}-1
\]
where $\sigma_{1}^{\mathrm{c}}$ the complement of $\sigma_{1}$ by the coefficient of $\alpha_{2}$.\\\\\\
ii) Roots set for for the sub-field $L_{2}$ that generated by the function $\sigma_{2}(z)$\\
\[{ }_{\sigma_{2}}^{L_{2,q'}} z_{s_{q'}}^{\mu}=\sigma_{2}^{-1}\left(\frac{\alpha_{0}}{\alpha_{2}}-\frac{1}{\alpha_{2}} \sigma_{2}^{c}\left(\sigma_{2}^{-1}\left(\frac{\alpha_{0}}{\alpha_{2}}-\frac{1}{\alpha_{2}} \sigma_{2}^{c}\left(\sigma_{2}^{-1}\left(\frac{\alpha_{0}}{\alpha_{2}}-\ldots- \right.\right.\right.\right.\right.\]
\[\left.\left.\left.\left.\left.-\frac{1}{\alpha_{2}} \sigma_{2}^{c}\left(\sigma_{2}^{-1}\left(\frac{\alpha_{0}}{\alpha_{2}}-\frac{1}{\alpha_{2}} \cdot \sum\limits_{\mathrm{i}=1}^{1} \alpha_{\mathrm{i}} \sigma_{\mathrm{i}}\left(\mathrm{z}_{0}\right)\right)\right)\right)\right)\right)\right)\right), k = 2, \ 0\leq s\leq \infty , q\in \mathbb{N}^{+}, \]
Apply: $\alpha_{2} \sigma_{2}(z)+\sigma_{2}^{c}(z)-\alpha_{0}=0 \wedge \sigma_{2}^{c}(z)=\sum\limits_{i=1}^{1} \alpha_{1} \sigma_{i}(z)=\alpha_{1} \cdot z^{r_{1}} \cdot b^{z^{q_{1}}}$
\[
\sigma_{2}^{-1}(z)=\left\{\left(\frac{r_{2}}{q_{2} \cdot \log b}\right) \cdot W\left(s'_{1}, \frac{\log b}{\frac{r_{2}}{q_{2}}} \cdot u^{\frac{q_{2}}{r_{2}}}\right)\right\}^{\frac{1}{q_{2}}} \cdot e^{\frac{2 s'_{2} \pi \cdot I}{q_{2}}},\left\{s_{2}, s_{1} \in Z\right\}, s_{2}=0 \div q_{2}-1
\]
where $\sigma_{2}^{\mathrm{c}}$ the complement of $\sigma_{2}$ by the coefficient of $\alpha_{1}$.\\\\
\textbf{IV.1.6. We will look at a real example.} We will look at a real example. Solve equation $x \cdot \alpha^{\frac{1}{x}}+\frac{1}{x} \alpha^{x}=2 a$.\\\\ 
According to the theory developed above we will have\\\\
i) Roots set for for the sub-field $L_{1}$ that generated by the function $\sigma_{1}(z)$\\\\
We have the constants: 
\[\alpha_{1}=1, \alpha_{2}=1, r_{1}=1, q_{1}=-1, a=\alpha, b=\alpha, \alpha_{0}=2 \alpha, r_{2}=-1, q_{2}=1 \]
\[{ }_{\sigma_{1}}^{L_{1,q}} z_{s_{q}}^{\mu}=\sigma_{1}^{-1}\left(\frac{\alpha_{0}}{\alpha_{1}}-\frac{1}{\alpha_{1}} \sigma_{1}^{c}\left(\sigma_{1}^{-1}\left(\frac{\alpha_{0}}{\alpha_{1}}-\frac{1}{\alpha_{1}} \sigma_{1}^{c}\left(\sigma_{1}^{-1}\left(\frac{\alpha_{0}}{\alpha_{1}}-\ldots- \right.\right.\right.\right.\right.\]
\[\left.\left.\left.\left.\left.-\frac{1}{\alpha_{1}} \sigma_{1}^{c}\left(\sigma_{1}^{-1}\left(\frac{\alpha_{0}}{\alpha_{1}}-\frac{1}{\alpha_{1}} \cdot \sum\limits_{\mathrm{i}=2}^{2} \alpha_{\mathrm{i}} \sigma_{\mathrm{i}}\left(\mathrm{z}_{0}\right)\right)\right)\right)\right)\right)\right)\right), 1 \leq k \leq 2, \ 0\leq s\leq \infty , q\in \mathbb{N}^{+},\]\\
Apply: $\alpha_{1} \sigma_{1}(z)+\sigma_{1}^{c}(z)-\alpha_{0}=0 \wedge \sigma_{1}^{c}(z)=\sum\limits_{i=2}^{2} \alpha_{2} \sigma_{i}(z)=\alpha_{2} \cdot z^{r_{2}} \cdot b^{z^{q_{2}}}$
\[\sigma_{1}^{-1}(z)=\left\{\left(\frac{r_{1}}{q_{1} \cdot \log a}\right) \cdot W\left(s_{1}, \frac{\log a}{\frac{r_{1}}{q_{1}}} \cdot u^{\frac{q_{1}}{r_{1}}}\right)\right\}^{\frac{1}{q_{1}}} \cdot e^{\frac{2 s_{2} \pi \cdot I}{q_{1}}},\left\{s_{2}, s_{1} \in Z\right\}, s_{2}=0 \div q_{1}-1\]
where $\sigma_{1}^{\mathrm{c}}$ the complement of $\sigma_{1}$ by the coefficient of $\alpha_{2}$.\\\\
ii) Roots set for for the sub-field $L_{2}$ that generated by the function $\sigma_{2}(z)$\\\\
We have the constants: 
\[\alpha_{1}=1, \alpha_{2}=1, r_{1}=1, q_{1}=-1, a=\alpha, b=\alpha, \alpha_{0}=2 \alpha, r_{2}=-1, q_{2}=1 \]
\[{ }_{\sigma_{2}}^{L_{2,q'}} z_{s_{q'}}^{\mu}=\sigma_{2}^{-1}\left(\frac{\alpha_{0}}{\alpha_{2}}-\frac{1}{\alpha_{2}} \sigma_{2}^{c}\left(\sigma_{2}^{-1}\left(\frac{\alpha_{0}}{\alpha_{2}}-\frac{1}{\alpha_{2}} \sigma_{2}^{c}\left(\sigma_{2}^{-1}\left(\frac{\alpha_{0}}{\alpha_{2}}-\ldots- \right.\right.\right.\right.\right.\]
\[\left.\left.\left.\left.\left.-\frac{1}{\alpha_{2}} \sigma_{2}^{c}\left(\sigma_{2}^{-1}\left(\frac{\alpha_{0}}{\alpha_{2}}-\frac{1}{\alpha_{2}} \cdot \sum\limits_{\mathrm{i}=1}^{1} \alpha_{\mathrm{i}} \sigma_{\mathrm{i}}\left(\mathrm{z}_{0}\right)\right)\right)\right)\right)\right)\right)\right), k = 2, \ 0\leq s\leq \infty , q\in \mathbb{N}^{+}, \]
Apply: $\alpha_{2} \sigma_{2}(z)+\sigma_{2}^{c}(z)-\alpha_{0}=0 \wedge \sigma_{2}^{c}(z)=\sum\limits_{i=1}^{1} \alpha_{1} \sigma_{i}(z)=\alpha_{1} \cdot z^{r_{1}} \cdot b^{z^{q_{1}}}$
\[
\sigma_{2}^{-1}(z)=\left\{\left(\frac{r_{2}}{q_{2} \cdot \log b}\right) \cdot W\left(s'_{1}, \frac{\log b}{\frac{r_{2}}{q_{2}}} \cdot u^{\frac{q_{2}}{r_{2}}}\right)\right\}^{\frac{1}{q_{2}}} \cdot e^{\frac{2 s'_{2} \pi \cdot I}{q_{2}}},\left\{s_{2}, s_{1} \in Z\right\}, s_{2}=0 \div q_{2}-1
\]
where $\sigma_{2}^{\mathrm{c}}$ the complement of $\sigma_{2}$ by the coefficient of $\alpha_{1}$.
\newpage\noindent
To get tangible results, we need to give a value for $\alpha$. So we randomly choose $\alpha=5$. Therefore with the theory, we will therefore construct 2 programs depending on each of the assemblies we are considering.\\\\
\textbf{I) Program 1.} Calculation of Roots set for the sub-field $L_{1}$ that generated by the function $\sigma_{1}(z)$\\\\
Constants: $\mathrm{a}:=5 ; \mathrm{n} 1:=1 ; 11:=-1 ; \mathrm{n} 2:=-1 ; 12:=1; \mathrm{b}:=\mathrm{a} ; \mathrm{a} 1:=1 ; \mathrm{k} 2:=1 ; \mathrm{a} 0:=2 \mathrm{a} ; \mathrm{a} 2:=1$;\\[2pt]
Clear [xr,y];Off[FindRoot::lstol,Divide::infy,Infinity::indet,FindRoot::jsing];\\[2pt]
$\mathrm{h} 2\left[\mathrm{x}_{-}\right]=\left(-\mathrm{a} 2^{*} \mathrm{x}^{\wedge} \mathrm{n} 2^{*} \mathrm{~b}^{\wedge}\left(\mathrm{x}^{\wedge}(12)\right)+\mathrm{a} 0\right) / \mathrm{a} 1$\\[2pt]
$\mathrm{h}\left[\mathrm{x}_{-}\right]:=\mathrm{a} 1^{*} \mathrm{x}^{\wedge} \mathrm{n} 1^{*} \mathrm{a}^{\wedge}\left(\mathrm{x}^{\wedge}(11)\right)+\mathrm{a} 2^{*} \mathrm{x} \wedge \mathrm{n} 2^{*} \mathrm{~b}^{\wedge}\left(\mathrm{x}^{\wedge}(12)\right)-\mathrm{a} 0$;\\[2pt]
For $[\mathrm{k} 1=-20, \mathrm{k} 1<=20, \mathrm{k} 1++$,\\[2pt]
f[u$_{-}$]:=\\[2pt]
$\mathrm{N}\left[\left(\text { ProductLog }\left[\mathrm{k} 1, \mathrm{n} 1 / 11 * \log [\mathrm{a}]^{*} \mathrm{u}^{\wedge}(\mathrm{n} 1 / 11)\right]^{*}((11 / \mathrm{n} 1) / \log [\mathrm{a}])\right)^{\wedge}(1 / \mathrm{n} 1)^{*} \operatorname{Exp}\left[2^{*} \mathrm{k} 2^{*} I^{*} \pi / \mathrm{n} 1\right]\right]$;\\[2pt]
$\mathrm{xr}=\mathrm{N}[\mathrm{Nest}[\mathrm{f}[\mathrm{h} 2[\#]] \&, 1 / 100,5]]$;\\
$\mathrm{FQ} 2=\mathrm{N}[\mathrm{xr}, 15]$\\[2pt]
Print $[" x(", k 1, ")=", s=N[y] / .$ FindRoot $[h[y]==0,\{y, F Q 2\}$,\\[2pt]
WorkingPrecision$\rightarrow$20],",","Approach $=$ ", $\mathrm{N}[\mathrm{h}[\mathrm{s}], 2]]]$\\\\
\textbf{Results 1:}\\\\
$\begin{aligned}
&x(-10)=4.56016142960097+39.0281378491641 \ \mathrm{I}\\
&x(-9)=4.43082289116071+35.1274487203240 \ \mathrm{I}\\
&x(-8)=4.28673006208215+31.2281347715740 \ \mathrm{I}\\
&x(-7)=4.12419318593592+27.3309211920972 \ \mathrm{I}\\
&x(-6)=3.93801587719292+23.4370540435161 \  \mathrm{I}\\
&x(-5)=3.72062664044545+19.5488002623928 \ \mathrm{I}\\
&x(-4)=3.46061196593289+15.6705560300157 \ \mathrm{I}\\
&x(-3)=3.14039657964064+11.8113468918753 \ \mathrm{I}\\
&x(-2)=2.73364128469454+7.98932994522154 \ \mathrm{I}\\
&x(-1)=2.20125675516346+4.22307532327606 \ \mathrm{I}\\
&x(0)=1.00000000000000\\
&x(1)=2.20125675516346-4.22307532327606 \  \mathrm{I}\\
&x(2)=2.73364128469454-7.98932994522154  \ \mathrm{I}\\
&x(3)=3.14039657964064-11.8113468918753 \  \mathrm{I}\\
&x(4)=3.46061196593289-15.6705560300157  \ \mathrm{I}\\
&x(5)=3.72062664044545-19.5488002623928  \ \mathrm{I}\\
&x(6)=3.93801587719292-23.4370540435161  \ \mathrm{I}\\
&x(7)=4.12419318593592-27.3309211920972 \ \mathrm{I}\\
&x(8)=4.28673006208215-31.2281347715740 \ \mathrm{I}\\
&x(9)=4.43082289116071-35.1274487203240 \ \mathrm{I}\\
&x(10)=4.56016142960097-39.0281378491641 \ \mathrm{I}
\end{aligned}$\\\\
20 complex roots (per 10 conjugates) in absolute ascending order with approximation $10^{\wedge}(-20)$.\\\\\\
\textbf{II. Program 2.} Calculation of Roots set for for the sub-field $L_{2}$ that generated by the function $\sigma_{2}(z)$\\\\
Constants: $\mathrm{a}:=5 ; \mathrm{n} 1:=1 ; 11:=-1 ; \mathrm{n} 2:=-1 ; 12:=1 ; \mathrm{b}:=\mathrm{a} ; \mathrm{a} 1:=1 ; \mathrm{k} 2:=0 ; \mathrm{a} 0:=2 \mathrm{a} ; \mathrm{a} 2:=1$;\\[2pt]
Clear[xr,y];Off[FindRoot::lstol,Divide::infy,Infinity::indet,FindRoot::jsing];\\[2pt]
$\mathrm{h} 2\left[\mathrm{x}_{-}\right]:=\left(-\mathrm{a} 1^{*} \mathrm{x}^{\wedge} \mathrm{n} 1^{*} \mathrm{a}^{\wedge}\left(\mathrm{x}^{\wedge}(11)\right)+\mathrm{a} 0\right) / \mathrm{a} 2$;\\[2pt]
$\mathrm{h}[\mathrm{x}]:=\mathrm{a} 1^{*} \mathrm{x}^{\wedge} \mathrm{n} 1^{*} \mathrm{a}^{\wedge}\left(\mathrm{x}^{\wedge}(11)\right)+\mathrm{a} 2^{*} \mathrm{x} \wedge \mathrm{n} 2^{*} \mathrm{~b}^{\wedge}\left(\mathrm{x}^{\wedge}(12)\right)-\mathrm{a} 0$;\\[2pt]
For[k1 $=-20, \mathrm{k} 1<=20, \mathrm{k} 1++$,\\[2pt]
$\mathrm{f}[\mathrm{u}]$ : =\\[2pt]
$\mathrm{N}\left[\left(\operatorname{ProductLog}\left[\mathrm{k} 1, \mathrm{n} 2 / 12^{*} \log [\mathrm{b}]^{*} \mathrm{u}^{\wedge}(\mathrm{n} 2 / 12)\right]^{*}((12 / \mathrm{n} 2) / \log [\mathrm{b}])\right)^{\wedge}(1 / \mathrm{n} 2)^{*} \operatorname{Exp}\left[2^{*} \mathrm{k} 2^{*} \mathrm{I}^{*} \pi / \mathrm{n} 2\right]\right]$;\\[2pt]
$\mathrm{xr}=\mathrm{N}[\mathrm{Nest}[\mathrm{f}[\mathrm{h} 2[\#]] \&, 100,7]]$\\[2pt]
$\mathrm{FQ} 2=\mathrm{N}[\mathrm{xr}, 15]$;\\[2pt]
Print $[" x(", k 1, ")=", s=N[y]/$.FindRoot$[h[y]==0,\{y, F Q 2\}$,\\[2pt]
WorkingPrecision$\to20], ", "$ "Approach $=$ ", $N[h[s], 2]]]$\\\\
\textbf{Results 2:}\\\\
$\begin{aligned}
&x(10)=0.00295349037853237+0.0252774449784375 \ \mathrm{I}\\
&x(9)=0.00353456403606972+0.0280219318116643 \ \mathrm{I}\\
&x(8)=0.00431446145318829+0.0314301534678917 \  \mathrm{I}\\
&x(7)=0.00539824183640218+0.0357740085279072 \  \mathrm{I}\\
&x(6)=0.00697236486206852+0.0414959454659245 \  \mathrm{I}\\
&x(5)=0.00939555374766368+0.0493658249852404 \  \mathrm{I}\\
&x(4)=0.0134370696673487+0.0608465656289249 \ \mathrm{I}\\
&x(3)=0.0210242783146740+0.0790744219809211 \ \mathrm{I}\\
&x(2)=0.0383388136452226+0.112048875481559  \ \mathrm{I}\\
&x(1)=1.00000005670143+3.53918967781277^{*} 10^{-9} \ \mathrm{I}\\
&x(0)=1.0000000000000\\
&x(-1)=1.00000005670143-3.53918967781277^{*} 10^{-9} \ \mathrm{I}\\
&x(-2)=0.0383388136452226-0.112048875481559 \ \mathrm{I}\\
&x(-3)=0.0210242783146740-0.0790744219809211 \ \mathrm{I}\\
&x(-4)=0.0134370696673487-0.0608465656289249 \  \mathrm{I}\\
&x(-5)=0.00939555374766368-0.0493658249852404 \ \mathrm{I}\\
&x(-6)=0.00697236486206852-0.0414959454659245 \mathrm{I}\\
&x(-7)=0.00539824183640218-0.0357740085279072 \ \mathrm{I}\\
&x(-8)=0.00431446145318829-0.0314301534678917  \ \mathrm{I}\\
&x(-9)=0.00353456403606972-0.0280219318116643\  \mathrm{I}\\
&x(-10)=0.00295349037853237-0.0252774449784375 \ \mathrm{I}
\end{aligned}$\\\\
Also 20 complex roots (per 10 conjugates) in absolute ascending order with approximation $10^{\wedge}(-20)$.\\\\
As it is known we have only 2 infinite categories per root case. By combining the values between the 2 approximations and the initial value, we achieve the best possible refinement of the roots. Usually we have 2 categories of initial values. Those that are in the interval $[0,1)$ and those that belong to the interval (1,infinity).\\\\\\
\textbf{IV.2. Solving of polynomial equation.}\\\\
For the determination of the roots of polynomial equation $f(x)=\sum\limits_{i=1}^{n} a_{i}\sigma_{i}(x)-a_{0}=0, \ \{x,a_{i},a_{0}\in\mathbb{C}, i\in\mathbb{N}^{+} \}$, which consists from distinct monomials functions $\sigma_{i}(x), 1 \leq i \leq n \wedge a_{i} \in C$, exist a repeating relation for every subfield of roots $G_{1}, G_{2}, \ldots, G_{r}, \ldots G_{n}$ and $\mathrm{n}$ the counting of the roots for any $G_{r}, r \in\{1,2, \ldots, n\}$. If we dissociate the initial function $f(x)=\sum\limits_{i=1}^{n} a_{i} \sigma_{i}(x)-a_{0}=0$ if we dissociate the initial function $f(x)=a_{n}\sigma_{n}(x)+\sum\limits_{i=1}^{n-1} a_{i} \sigma_{i}(x)-a_{0}=0$ to 2 parts then we have the form ${ }_{r} x_{j+1}^{r}=-\frac{1}{a_{r}} \sum\limits_{i \neq r, i=1}^{n} a_{i} x_{j+1}^{i}+\dfrac{a_{0}}{a_{r}} \Rightarrow {}_{r}x_{j+1}=e^{2 k \pi i / r} \cdot\left(-\dfrac{1}{a_{r}} \sum\limits_{i \neq r, i=1}^{n} a_{i} x_{j}^{i} + \dfrac{a_{0}}{a_{r}}\right)^{1 / r}, \ \{ a_{r}\in\mathbb{C},k=0,1,2,..n-1,1\leq r\leq n\in\mathbb{N}^{+}\}$ according to the iterative procedure. The process starts with $\mathrm{r}=\mathrm{n}$, field $\left\{G_{n}\right\}$ and continues with $\mathrm{r}=\mathrm{n}-1, \mathrm{r}=\mathrm{n}-2$ fields $G_{n-1}, G_{n-2}...$ etc until they arrive all the roots of equation and the field G.\\\\
\textbf{Example}\\\\
Calculate the roots of the equation $a_{n}x^{n}-a_{1}x-a_{0}=0, \{ a_{n}=1,a_{1},a_{0}\in\mathbb{C}\}$ with the method GR$_{\text{T}}$E.\\\\
According to the iterative procedure, apply:
\[\left({ }_{n} x_{j+1}\right)^{n}=\frac{1}{a_{n}} \sum_{i =1}^{1} a_{i} x_{j+1}^{i}+\dfrac{a_{0}}{a_{n}} \Rightarrow {}_{r}x_{j+1}=e^{2 k \pi i / n} \cdot\left(\frac{1}{a_{n}} \sum_{i =1}^{1} a_{i} x_{j}^{i}+\dfrac{a_{0}}{a_{n}}\right)^{1 / n},\]
\[\{a_{n}=1,a_{1},a_{0}\in\mathbb{C},k=0,1,2,..n-1,n\in\mathbb{N}^{+} \}\]
with method GR$_{\text{T}}$E a simple program will be ,if initial value is $x 0$ and for 4 repeated times:\\\\
$\begin{aligned}
&\mathrm{fl}\left[\mathrm{x}_{-}\right]:=\mathrm{x}^{\wedge} \mathrm{n} ; \\[2pt]
&\mathrm{f} 2\left[\mathrm{x}_{-}\right]=\mathrm{a} 1^{*} \mathrm{x}+\mathrm{a} 0 ; \\[2pt]
&\mathrm{h} 2\left[\mathrm{x}_{-}\right]:=\mathrm{f} 2[\mathrm{x}] ; \\[2pt]
&\mathrm{h}[\mathrm{x}]:=\mathrm{fl}[\mathrm{x}]-\mathrm{f} 2[\mathrm{x}] ; \\[2pt]
&\mathrm{f}\left[\mathrm{u}_{-}\right]:=\mathrm{u} \wedge(1 / \mathrm{n})^{*} \operatorname{Exp}\left[2 * \mathrm{k}^{*} \pi^{*} \mathrm{I} / \mathrm{n}\right] ; \\[2pt]
&\mathrm{xn}=\mathrm{Nest}[\mathrm{f}[\mathrm{h} 2[\#]] \&, \mathrm{x} 0,4]
\end{aligned}$\\\\
\\The answer in general form is:
\[
x_{4}=e^{2 \pi k I / n}\left(a_{0}+a_{1} e^{2 \pi k I / n}\left(a_{0}+a_{1} e^{2 \pi k I / n}\left(a_{0}+a_{1} e^{2 \pi k I/ n}\left(a_{0}+a_{1} x_{0}\right)^{1 / n}\right)^{1 / n}\right)^{1 / n}\right)^{1 / n} \dots, 
\]
\[\text { with } \mathrm{k}=0,1, \ldots \mathrm{n}-1 
\]
The method can, of course, be combined \textbf{with the Newton method} because it is local, yielding results with a rather very good approach. Complete examples of such can be seen at the end of the work.\\\\\\
\textbf{IV.3. Equations of category with Transcendental or Irrational roots to polynomial form.}\\\\
We look at the form $x^{m}-p x^{s}+q=0,\{(m, s) \in R,(p, q) \in C$) (2) again in case c, but now with an extension of the set $Z$ to $R$ of G.R.P-N. The process is similar to $b$, we just have sub-cases that depend with the number of roots in relation to the maximum exponent $\mathrm{m}$ of the equation. We will see the whole process with examples in particular. The difference is to the fact that the number of roots in the transcendental is the integer part of the Real maximum exponent with a value $\pm 1$, the final value of which depends on the \textbf{sign of the last fixed term} $\mathbf{q}$, of the equation. So if it is a sign of $\mathbf{q}<\mathbf{0}$ then the count of roots $x$ of transcendental equation (2) will be $x=\operatorname{IntegerPart}[\max \{m, s\}]+1$, different if $\mathbf{q}>\mathbf{0}$ then $x=$ IntegerPart $[\max \{m, s\}]$. This fact appears to be clear in all cases where we can meet. And because we are looking always the binomial $\mid x^{\alpha}=u$, if $a=\max \{m, s\}$ that like if $x^{\alpha}=u \Rightarrow x=e^{\frac{2 \cdot k \pi \cdot I}{\alpha}} u^{1 / \alpha}, u \in C$ we distinguish two basic cases:\\
\begin{enumerate}[label=\Roman*)]
    \item  $[a]=$ IntegerPart$[a]=2 n+1, n \in Z$
    \begin{enumerate}[label=\roman*)]
        \item $\mathbf{u>0}$ then $k=0, \pm 1, \ldots \pm([a]-1) / 2$ and
        \item $\mathbf{u}<\mathbf{0}$ then $k=0, \pm 1, \ldots, \pm([a]-1) / 2,-([a]+1) / 2$
    \end{enumerate}
    \item $[a]=$ IntegerPant $[a]=2 n, n \in Z$
    \begin{enumerate}[label=\roman*)]
        \item[iii)] $\mathbf{u}>0$ then $k=0, \pm 1, \ldots \pm([a]) / 2$ and iv) $\mathbf{u}<\mathbf{0}$ then $k=0, \pm 1, \ldots,-([a]) / 2$
    \end{enumerate}
\end{enumerate}
\vspace{\baselineskip}
\textbf{Examples}\\\\
I) $[a]=$ IntegerPart$[a]=2 n+1, n \in Z$\\\\
i) $u>0$\\\\
The binomial $z^{\pi}=2$ it has complex roots. This equation is a transcendental one and accepts the following solutions: Apply $[\pi]=3$ therefore $k=0, \pm 1$ because $k \leq([a]-1) / 2=1$. These solutions are given more analytically by the relations below:
\[
\begin{aligned}
&z_{0}=2^{1 / \pi} e^{0}=1.24686 \\
&z_{1}=2^{1 / \pi} e^{2 \pi i / \pi}=-0.51888+1.33774 i \\
&z_{-1}=2^{1 / \pi} e^{-2 \pi i / \pi}=-0.51888-1.33774 i
\end{aligned}
\]\\
ii) $u<0$\\\\
The binomial $z^{\pi}=-2$ it has only complex roots. This equation has solutions: Because apply $[\pi]=3$ therefore $k=0, \pm 1,-2$. We accept more analytically by the relations below:
\[
\begin{aligned}
&z_{0}=(-2)^{1 / \pi} e^{0}=0.67368+1.04920 i, \ \  z_{-1}=(-2)^{1 / \pi} e^{-2 \pi i / \pi}=0.67368+1.04920 i\\[8pt]
&z_{1}=(-2)^{1 / \pi} e^{2 \pi i / \pi}=-1.23439+0.17595 i, \ \  z_{-2}=(-2)^{1 / \pi} e^{-2.2 \pi i / \pi}=-1.23439-0.17595 i
\end{aligned}
\]
We have 4 solutions\\\\
II) $[a]=$ IntegerPart$[a]=2 n, n \in Z$\\\\
i) $u>0$\\\\
The binomial $z^{\pi+1}=2$ we ask complex roots. We accept the following solutions: Apply $[\pi+1]=4$ therefore $k=0, \pm 1, \pm 2$ because $k \leq([a]) / 2=2$. These solutions are given more analytically by the relations below:
\[
\begin{aligned}
&z_{0}=2^{1 /(\pi+1)} e^{0}=1.18218 \\
&z_{1,-1}=2^{1 /(\pi+1)} e^{\pm 2 \pi i /(\pi+1)}=0.06345546 \pm 1.18047843 \mathrm{i}, \ \  z_{2,-2}=2^{1 / \pi} e^{\pm 4 \pi i / \pi}=-1.17537055 \pm 0.12672797\mathrm{i}
\end{aligned}
\]
We have 5 solutions.\\\\
ii) $u<0$\\\\
The binomial $z^{\pi+1}=-2$ it has only complex roots. We get the following solutions: Apply $[\pi+1]=4$ therefore $k=0, \pm 1,-2$. These solutions are given more analytically:
\[\begin{aligned}
&z_{0}=2^{1 /(\pi+1)} e^{0}=0.85807105+0.81318508 \mathrm{i} \\
&z_{-1}=2^{1 /(\pi+1)} e^{-2 \pi i /(\pi+1)}=0.85807105+0.81318508 \mathrm{i}, \\
&z_{1}=2^{1 / \pi} e^{2 \pi i / \pi}=-0.76595450+0.90048299 \mathrm{i} \\
&z_{-2}=2^{1 / \pi} e^{-4 \pi i / \pi}=-0.76595450-0.90048299 \mathrm{i}
\end{aligned}\]
Therefore we have 4 solutions.\\
At the same time, we also have the solution of the principal transcendental equation we shall need subsequently for solving the more general forms of transcendental equations.\\\\
\textbf{IV.4. Equation category to full trinomial} of the form $a^{x}+p \cdot b^{x}-q=0,\{a \neq 1\wedge  b \neq 1,a, b, p, q \in \mathbb{C}\}$.\\\\
We solve the two transcendental binomials $a^{x}=u$ and $b^{x}=u$ that arise if we solve each term individually. The resulting fields for each binomial therefore arise are.. $x=(\log u+2 \cdot k \cdot \pi \cdot I) / \log a$ and $x=(\log u+2 \cdot k \cdot \pi \cdot I) / \log b .$\\\\
We have 2 Programs in mathematica, for equation of form $a^{x}+p \cdot b^{x}-q=0,\{a \neq 1\wedge b \neq 1,a, b, p, q \in \mathbb{C}\}$ of G.R.P-N Method. Will be:\\\\
\textbf{1).}\\\\
Constants: a: $=7 / 5 ; b:=4 / 5 ; v:=20 ; k:=0,+1,-1+2,-2 \ldots$ infinity\\[2pt]
$\mathrm{p}:=1 ; \mathrm{q}:=1$;\\[2pt]
h1[x$_{-}$]:=(a) $x+p^{*}(b)^{\wedge} x-q$\\[2pt]
h2[x$_{-}$]:=$-p^{*}(b)^{\wedge} x+q$\\[2pt]
$\mathrm{k}:=1 ; \mathrm{f}\left[\mathrm{u}_{-}\right]:=\left(\text{Log} [\mathrm{u}]+\left(2^{*} \mathrm{k}^{*} \pi^{*} \mathrm{I}\right)\right) / \text{Log} [\mathrm{a}]$\\[2pt]
$\mathrm{xr}=\mathrm{Nest}[\mathrm{f}[(\mathrm{h} 2[\#])] \&, 1, \mathrm{v}]$;\\[2pt]
$\mathrm{s} 2=\mathrm{N}[\mathrm{Abs}[\mathrm{xr}], 5]$\\[2pt]
$\mathrm{FQ} 2=\mathrm{N}[\mathrm{xr}, 10]$\\[2pt]
Print["Approach=",f1[d]/.FindRoot[ $[\mathrm{f1}[\mathrm{d}]=0,\{\mathrm{d}, \mathrm{fq}\}$,\\[2pt]
WoringPrecision$\rightarrow 10], ", x(", k, ")=", N[d] /$ FindRoot $[\mathrm{f1}[\mathrm{d}]=0$,\\[2pt]
$\{\mathrm{d}, \mathrm{fq}\}\}$, WorkingPrecision$\rightarrow 10]] / /$ Timing\\\\
According to the program and these values for $k$ i get:\\\\
\textbf{Results 3:}\\\\
$\begin{aligned}
&\mathrm{k}=0, \mathrm{x}=-0.444738772955+3.753482934416 \mathrm{i} \\
&\mathrm{k}=1, \mathrm{x}=1.3939724221+17.6753885427 \mathrm{i} \\
&\mathrm{k}=-1, \mathrm{x}=1.3939724221-17.6753885427 \mathrm{i} \\
&\mathrm{k}=2, \mathrm{x}=1.3685154337+38.4113462022 \mathrm{i} \\
&\mathrm{k}=-2, \mathrm{x}=1.3685154337-38.4113462022 \mathrm{i}
\end{aligned}$\\\\
\hspace{7cm} Infinity Roots\\\\
Significance has all the roots and which, as we know, is infinite. To complete the set of roots we should look at the other group of roots corresponding to the other field.
\newpage\noindent
\textbf{2).}\\\\
Constants: $a:=7 / 5 ; b:=4 / 5 ; v:=20 ; k:=0,+1,-1+2,-2$...infinity, $p:=1 ; q:=1$;\\[2pt]
$\mathrm{h} 1\left[\mathrm{x}_{-}\right]:=(\mathrm{a})^{\wedge} \mathrm{x}+\mathrm{p}^{*}(\mathrm{~b})^{\wedge} \mathrm{x}-\mathrm{q}$\\[2pt]
h$2\left[\mathrm{x}_{-}\right]:=-1 / p^{*}(a)^{\wedge} x+q / p$\\[2pt]
$\mathrm{k}:=1 ; \mathrm{f}\left[\mathrm{u}_{-}\right]:=\left(\text{Log} [\mathrm{u}]+\left(2^{*} \mathrm{k}^{*} \mathrm{\pi}^{*} \mathrm{I}\right)\right) / \text{Log} [\mathrm{b}$\\[2pt]
$\mathrm{xr}=$ Nest $[\mathrm{f}[(\mathrm{h} 2[\#])] \&, 1, \mathrm{v}]$;\\[2pt]
$\mathrm{s} 2=\mathrm{N}[\mathrm{Abs}[\mathrm{xr}], 5]$\\[2pt]
$\mathrm{FQ} 2=\mathrm{N}[\mathrm{xr}, 12]$\\[2pt]
Print["Approach=",f1[d]/.FindRoot $[\mathrm{fl}[\mathrm{d}]=0,\{\mathrm{d}, \mathrm{fq}\}$,\\[2pt]
WoringPrecision$\rightarrow 12], ", \mathrm{x}(", \mathrm{k}, ")=", \mathrm{~N}[\mathrm{~d}] /$.FindRoot $[\mathrm{fl}[\mathrm{d}]=0$,\\[2pt]
$\{\mathrm{d}, \mathrm{fq}\}$, WorkingPrecision$\rightarrow$12]]//Timing\\\\
According to the program and these values for $k$ i get:\\\\
\textbf{Results 4:}\\\\
$\begin{aligned}
&\mathrm{k}=0, \mathrm{x}=-1.89895932741-28.1072354396 \mathrm{i} \\
&\mathrm{k}=1, \mathrm{x}=-1.89895932741+28.1072354396 \mathrm{i} \\
&\mathrm{k}=2, \mathrm{x}=-0.528066150559-59.9081922009 \\
&\mathrm{k}=-2, \mathrm{x}=-0.528066150559-59.9081922009 \\
\end{aligned}$\\\\
Infinity Roots\\\\
In roots calculation many times one of the two cases does not yield roots. Generally, however, the whole of the roots are the union of the roots of the two fields.\\\\\\
\textbf{IV.5. Solving the equation $\boldsymbol{e^{x^{3}}-x^{2}-5=0, x \in C}$ in the set of complexes. [7,11]}\\\\
According to the method of periodic Radicals for the solution we are looking for:\\\\
\textbf{1). To determine the functional factors}\\\\
In our case we have 2 functional factors $\sigma_{1}(x)=e^{x^{3}} \wedge \sigma_{2}(x)=x^{2}$. Therefore the subfields that produced are $G_{1}, G_{2}$ where $G=G_{1} \cup G_{2}$ is the set of solutions of Roots of equation.\\\\
\textbf{2). To find each time the inverse function as to $\mathbf{x}$, functional factor $\boldsymbol{\sigma_{k}(x), 1 \leq k \leq 2}$}\\\\
For each functional factor resulting from the relation $\sigma_{k}(x)=u_{k} \Rightarrow x_{j}^{\sigma_{k}}=\sigma_{k}^{-1}\left(u_{k}\right)$ we take analytically:\\\\
$a) \therefore \sigma_{1}(x)=u_{1} \Rightarrow{ }_{\sigma_{1}} x_{j}=\sigma_{1}^{-1}\left(u_{1}\right) \Rightarrow e^{x^{3}}=u_{1} \Rightarrow{ }_{\sigma_{1}}^{m_{1}} x_{j}=(\log (u)+2 \cdot w \cdot \pi \cdot i)^{1 / 3} \cdot e^{\frac{2 m 1 \cdot \pi \cdot i}{3}}$, $0 \leq m_{1} \leq 2, w \in Z$\\\\
b) $\therefore \sigma_{2}(x)=u_{2} \Rightarrow \sigma_{\sigma_{2}} x_{j}=\sigma_{2}^{-1}\left(u_{2}\right) \Rightarrow \sigma_{2}^{m_{2}} x_{j}^{2}=u_{2} \Rightarrow \sigma_{\sigma_{2}}^{m_{2}} x_{j}=\sqrt{u_{2}} \cdot e^{\frac{2 m \cdot \pi \cdot i}{2}}, 0 \leq m_{2} \leq 1$\\\\
For both cases $\mathrm{j}$ denotes a repetitive process.\\\\
3. We determine approximate now in agreement with form ${ }_{\sigma_{k}}^{m_{k}} x_{\# +1}^{\sigma_{k}}=\Phi_{k, u}\left(\frac{1}{a_{k}}\left(-\sum\limits_{i=1, i \neq k}^{n} a_{i}\sigma_{i}(\# )-a_{0}\right)\right)$ where the symbol \# shows repetitive (periodic) procedure of values and $1<=k<=2$ and $\mu$ the number of repeats.
In our case we will have for each case separately:
\begin{enumerate}
    \item[a)] $\therefore$ Case:
    \begin{enumerate}
        \item[i)] $\therefore \Phi_{1, \mu}\left(u_{1}\right)=\left(\log \left(u_{1}\right)+2 \cdot w \cdot \pi \cdot i\right)^{1 / 3} \cdot e^{\frac{2 \cdot m_{1} \cdot \pi \cdot i}{3}}, 0 \leq m_{1} \leq 2, k=1, w \in Z$.
        \item[ii)]  $\therefore{ }_{\sigma_{1}}^{m_{1}} x_{\#+1}^{\mu}=\Phi_{1, \mu}\left(\frac{1}{a_{k}}\left(-\sum\limits_{i>1}^{2} a_{i} \sigma_{i}(\#)-a_{0}\right)\right)=\Phi_{1, \mu}\left(\frac{1}{a_{1}}\left(-\alpha_{2} \sigma_{2}(\#)-a_{0}\right)\right)=\Phi_{1, \mu}\left(\#^{2}+5\right)$.
    \end{enumerate}
    \item $\therefore$ Case: 
    \begin{enumerate}
        \item[i)] $\therefore \Phi_{2, \mu}\left(u_{2}\right)=\sqrt{u_{2}} \cdot e^{\frac{2 \cdot m_{2} \cdot \pi \cdot i}{2}}, 0 \leq m_{2} \leq 1, k=2$.
        \item[ii)] $\therefore{ }_{\sigma_{2}}^{m_{2}} x_{\#+1}^{\mu}=\Phi_{k, \mu}\left(\frac{1}{a_{k}}\left(-\sum\limits_{i=1, i \neq 2}^{2} a_{i} \sigma_{i}(\#)-a_{0}\right)\right)=\Phi_{2, \mu}\left(\frac{1}{a_{2}}\left(-\alpha_{1} \sigma_{1}(\#)-a_{0}\right)\right)=\Phi_{2, \mu}\left(e^{\#}-5\right)$.
    \end{enumerate}
\end{enumerate}
\vspace{\baselineskip}
\textbf{4. Numerical calculations for Complex roots:}\\\\
\textbf{4.1. Categories of complex roots corresponding to the case a:}\\\\
According to (case a), there are 3 categories of roots and more specifically we take.\\\\ \textbf{The basics Programs in mathematica}.\\\\
$\mathrm{d}:=6 ; \mathrm{w}:=0,1,2 ;\\[2pt]
\mathrm{h} 2[\mathrm{x}]:=\mathrm{x}^{\wedge} 2+5 \text {; }\\[2pt]
\mathrm{h}\left[\mathrm{x}_{-}\right]:=\mathrm{e}^{\mathrm{x}^{\wedge} 3}-\mathrm{x}^{\wedge} 2-5 \text {; }\\[2pt]
\text {For }[\mathrm{k}=-\mathrm{d}, \mathrm{k}<=\mathrm{d}, \mathrm{k}++\text {, }\\[2pt]
\mathrm{f}\left[\mathrm{u}_{-}\right]:=(2 \mathrm{I} \pi \mathrm{k}+\log {[\mathrm{u}])^{1/3}}*\operatorname{Exp}\left[2^{*}\mathrm{I}^{*}\pi^{*}\mathrm{w}/3\right];\\[2pt]
\mathrm{xr}=\mathrm{Nest}[\mathrm{f}[\mathrm{h} 2[\#]] \&, 1,15];\\[2pt]
\mathrm{FQ} 2=\mathrm{N}[\mathrm{xr}, 25];\\[2pt]
\operatorname{Print}[" \mathrm{x}(", \mathrm{k}, ")=", \mathrm{~N}[\mathrm{y}] / \text { FindRoot }[\mathrm{h}[\mathrm{y}]==0,\{\mathrm{y}, \mathrm{FQ} 2\} \text {, WorkingPrecision}\to20]]]$\\\\
\textbf{i. If $\boldsymbol{w=1}$ for the first 13 roots (from infinite complex Roots) we get.}\\\\
\textbf{Results 5:}\\\\
$\begin{aligned}
&x(-6)=-0.0541597908772764481+3.4430687795713415686 \text{ }\mathrm{I} \\
&\mathrm{x}(-5)=-0.0542272093684110848+3.2561465116970043623 \text{ } \mathrm{I} \\
&\mathrm{x}(-4)=-0.0522799127977140002+3.0446951231985700472 \text{ } \mathrm{I} \\
&\mathrm{x}(-3)=-0.044446198208043990301+2.7986415748311100303 \text{ } \mathrm{I} \\
&\mathrm{x}(-2)=0.073991657283343344493+2.2314465835450111411 \text{ } \mathrm{I} \\
&\mathrm{x}(-1)=-0.043290439396162264200+1.8563014885856711080 \text{ } \mathrm{I} \\
&\mathrm{x}(0)=-0.51628371223967864689+1.0193275613354981293 \text{ } \mathrm{I} \\
&\mathrm{x}(1)=-1.4767721464906745539+1.0770106461916786832 \text{ } \mathrm{I} \\
&\mathrm{x}(2)=-1.9196926480095895878+1.2656352955061933943 \text{ } \mathrm{I} \\
&\mathrm{x}(3)=-2.2255896653786903526+1.4135307520156373273\text{ } \mathrm{I} \\
&\mathrm{x}(4)=-2.4668779334538722208+1.5359648431918860639 \text{ }\mathrm{I} \\
&\mathrm{x}(5)=-2.6694198069847862912+1.6414453768855680311\text{ } \mathrm{I} \\
&\mathrm{x}(6)=-2.8457158247398611933+1.7347795001319450124\text{ } \mathrm{I}
\end{aligned}$\\\\\\
\textbf{ii.If $w=0$ for the first 13 roots (from infinite complex Roots) we get.}\\\\
\textbf{Results 6:}\\\\
$\begin{aligned}
&\mathrm{x}(-6)=2.9629015049765169380-1.6197988551044190385\text{ } \mathrm{I} \\ &\mathrm{x}(-5)=2.7970790204948353259-1.5157734611427453324 \text{ }\mathrm{I} \\ &\mathrm{x}(-4)=2.6084406987474839689-1.3956240678529720229 \text{ }\mathrm{I} \\ &\mathrm{x}(-3)=2.3870122571990767896-1.2512153777260514123 \text{ }\mathrm{I} \\ &\mathrm{x}(-2)=2.1133652182429052922-1.0649772354315062827 \text{ }\mathrm{I} \\
&\mathrm{x}(-1)=1.7402048569203715448-0.7832829915835906632\text{ }\mathrm{I} \\
\end{aligned}$\\
$\begin{aligned}
&\mathrm{x}(0)=1.2330961898317350649 \\
&\mathrm{x}(1)=1.7402048569203715448+0.7832829915835906632\text{ } \mathrm{I} \\ &\mathrm{x}(2)=2.1133652182429052922+1.0649772354315062827 \text{ }\mathrm{I} \\ &\mathrm{x}(3)=2.3870122571990767896+1.2512153777260514123\text{ } \mathrm{I} \\ &\mathrm{x}(4)=2.6084406987474839689+1.3956240678529720229 \text{ }\mathrm{I} \\ &\mathrm{x}(5)=2.7970790204948353259+1.5157734611427453324 \text{ }\mathrm{I} \\ &\mathrm{x}(6)=2.9629015049765169380+1.6197988551044190385 \text{ }\mathrm{I} \\
\end{aligned}$\\\\\\
To this category belongs the \textbf{only Real root} of the equation.\\\\
$x(0)=1.2330961898317350649$\\\\\\
\textbf{iii. If $\boldsymbol{w=2}$ for the first 13 (from infinite complex Roots) we get.}\\\\
\textbf{Results 7:}\\\\
$\begin{aligned}
&\mathrm{x}(-6)=-2.8457158247398611933-1.7347795001319450124\text{ } \mathrm{I} \\
&\mathrm{x}(-5)=-2.6694198069847862912-1.6414453768855680311 \text{ }\mathrm{I} \\
&\mathrm{x}(-4)=-2.4668779334538722208-1.5359648431918860639 \text{ }\mathrm{I} \\
&\mathrm{x}(-3)=-2.2255896653786903526-1.4135307520156373273\text{ } \mathrm{I} \\
&\mathrm{x}(-2)=-1.9196926480095895878-1.2656352955061933943 \text{ }\mathrm{I} \\
&\mathrm{x}(-1)=-1.4767721464906745539-1.0770106461916786832 \text{ }\mathrm{I} \\
&\mathrm{x}(0)=-0.51628371223967864689-1.0193275613354981293 \text{ }\mathrm{I} \\
&\mathrm{x}(1)=-0.043290439396162264200-1.8563014885856711080 \text{ }\mathrm{I} \\
&\mathrm{x}(2)=0.073991657283343344493-2.2314465835450111411 \text{ }\mathrm{I} \\
&\mathrm{x}(3)=-0.044446198208043990301-2.7986415748311100303 \text{ }\mathrm{I} \\
&\mathrm{x}(4)=-0.0522799127977140002-3.0446951231985700472 \text{ }\mathrm{I} \\
&\mathrm{x}(5)=-0.0542272093684110848-3.2561465116970043623 \text{ }\mathrm{I} \\
&\mathrm{x}(6)=-0.0541597908772764481-3.4430687795713415686 \text{ }\mathrm{I}
\end{aligned}$\\\\\\
\textbf{4.2.Categories of complex roots corresponding to the case b:}\\\\
According to (case b), there is 2 Categories of roots and namely we have.\\\\
The \textbf{basics Programs in mathematica}.\\\\
$\begin{aligned}
&\mathrm{d}:=0 ; \mathrm{h} 2\left[\mathrm{x}_{-}\right]:=\mathrm{e}^{\mathrm{x}^{\wedge} 3}-5 ; \mathrm{h}\left[\mathrm{x}_{-}\right]:=\mathrm{e}^{\mathrm{x}^{\wedge} 3}-\mathrm{x}^{\wedge} 2-5 \text {; }\\[2pt]
&\text {For }[\mathrm{k}=-\mathrm{d}, \mathrm{k}<=\mathrm{d}, \mathrm{k}++\text {, }\\[2pt]
&\mathrm{f}\left[\mathrm{u}_{-}\right]:=\mathrm{u}^{\wedge}(1 / 2)^{*} \mathrm{e}^{\wedge}(2 \mathrm{w} \pi \mathrm{i} / 2) ; \mathrm{xr}=\mathrm{Nest}\left[\mathrm{f}[\mathrm{h} 2[\#]] \&, \mathrm{x}_{0}, 6\right] ; \mathrm{FQ} 2=\mathrm{N}[\mathrm{xr}, 25] \text {; }\\[2pt]
&\text { Print }[" A p p r o a c h=", N[h[y], 30] / \text { FindRoot }[h[y]==0,\{y, F Q 2\} \text {, WorkingPrecision}\to 30] \text {, }\\[2pt]
&" \mathrm{x}(", \mathrm{k}, ")=", \mathrm{~N}[\mathrm{y}] / . \text { FindRoot }[\mathrm{h}[\mathrm{y}]==0,\{\mathrm{y}, \mathrm{FQ} 2\} \text {, WorkingPrecision}\to 30]]
\end{aligned}$\\\\\\
If $w=0 \& w=1$ we take 2 roots i.e we get 2 Complex Roots\\\\
Results:\\\\
$\begin{aligned}
&\left.x(0)=0.073991657283343344493-2.2314465835450111411 \text { I (with initial value } x_{0}=1\right) \\[2pt]
&\left.\bar{x}(0)=0.073991657283343344493+2.2314465835450111411 \text { I (with initial value } x_{0}=1 / 2\right)
\end{aligned}$\\\\\\
The Approach in both cases is $10^{-27} \pm 10^{-27} i$.\\\\\\
\textbf{IV.6. Solving the equation} $a^{x}+p \cdot b^{x}=q \cdot c^{x},\{x, a, b, c\} \in C$, \textbf{in the set of complexes.} In general, we have 1 methods \textbf{of solutions and 2 cases.} Both of these cases are generated using a logarithm. The functions that represent this method are cases that arise analytically per case:\\\\
We define the function:\\
\[f(x)=\left(\frac{a}{c}\right)^{x}+p \cdot\left(\frac{b}{c}\right)^{x}-q, \text{ } x \in C,\{a, b, c, p, q \in C \wedge a, b, c \neq 0 \}\]\\
a) Case: for $f_{1}(x)=q-p \cdot\left(\frac{b}{c}\right)^{x}$\\\\
i) $\therefore \Phi_{1, \mu}\left(u_{1}\right)=\frac{\log u_{1}+2 \cdot w \cdot \pi \cdot i}{\log \left(\frac{a}{c}\right)}, w \in Z$.\\\\[4pt]
b) Case: for $k=2, f_{2}(x)=\frac{q}{p}-\frac{1}{p} \cdot\left(\frac{a}{c}\right)^{x}$.\\\\
i) $\therefore \Phi_{2, \mu}\left(u_{2}\right)=\frac{\log u_{2}+2 \cdot w \cdot \pi \cdot i}{\log \left(\frac{b}{c}\right)}, w \in Z$.\\\\
ii) $\therefore \sigma_{2}^{m_{2}} x_{m+1}^{\mu}=\Phi_{k, \mu}\left(\frac{1}{a_{k}}\left(-\sum_{i \neq 2, i \geq \mathrm{I}}^{2} a_{i} \sigma_{i}(\#)-a_{0}\right)\right)=\Phi_{2, \mu}\left(f_{2}(\#)\right)$.\\\\\\
\textbf{Numerical calculations for Complex roots:}\\\\
\textbf{6.1. Categories of complex roots corresponding to the case a:}\\\\
\textbf{Constants:} 
$m_{1}=\frac{a}{c}=\sqrt[5]{43} ; m_{2}=\frac{b}{c}=-\sqrt[3]{31}$;\\[2pt]
$\mathrm{d}:=15 ; \mathrm{p}:=\pi ; \mathrm{q}:=\mathrm{e}$;\\[2pt]
$\mathrm{h} 1\left[\mathrm{x}_{-}\right]:=(\mathrm{m} 1)^{\wedge} \mathrm{x}+\mathrm{p}^{*}(\mathrm{m} 2)^{\wedge} \mathrm{x}-\mathrm{q};$\\[2pt]
$\mathrm{h} 2\left[\mathrm{x}_{-}\right]:=-\mathrm{p}^{*}(\mathrm{m} 2)^{\wedge} \mathrm{x}+\mathrm{q}$;\\[2pt]
For $[\mathrm{k}=0, \mathrm{k}<=\mathrm{d}, \mathrm{k}++$,\\[2pt]
$\mathrm{f}\left[\mathrm{u}_{-}\right]:=\left(\log [\mathrm{u}]+\left(2^{*} \mathrm{k}^{*} \pi^{*} \mathrm{I}\right)\right) / \log [\mathrm{m} 1]$;\\[2pt]
$\mathrm{xr}=$ Fold $[\mathrm{f}[\mathrm{h} 2[\#])] \&,-2$, Range[10]];\\[2pt]
FQ2=N $[x r, 10]$;\\[2pt]
Print["Approach=", $\mathrm{N}[\mathrm{h} 1[\mathrm{y}], 30] /$. FindRoot $[\mathrm{h} 1[\mathrm{y}]==0,\{\mathrm{y}, \mathrm{FQ} 2\}$, WorkingPrecision$\to 30], " \mathrm{x}(", \mathrm{k}, ")=", \mathrm{~N}[\mathrm{z}] /.$ FindRoot $[\mathrm{h} 1[\mathrm{z}]==0,\{\mathrm{z}, \mathrm{FQ} 2\}$, WorkingPrecision$\to 30]]]$\\\\
\textbf{For the first 15 complex roots (this case don't have spouses complex roots) we get.}\\\\
\textbf{Results 8:}\\\\
$\begin{aligned}
&x(0)=58.0487507536351199035328205089+7.61535990029372317354610795877 \text{ } \mathrm{I} \\
&\mathrm{x}(1)=1.32936292836325648633528741447+8.35263361944434710891163214802 \text{ } \mathrm{I} \\
&\mathrm{x}(2)=1.32936292837426009323673136176+16.7052672389407600393663302459 \text{ } \mathrm{I} \\
&\mathrm{x}(3)=1.32936292837426009323670176765+25.0579008584111400590493309734 \text{ } \mathrm{I} \\
&\mathrm{x}(4)=1.32936292837426009323670176765+33.4105344778815200787324412978 \text{ } \mathrm{I} \\
&\mathrm{x}(5)=1.32936292837426009323670176765+41.7631680973519000984155516223 \text{ } \mathrm{I} \\
&\mathrm{x}(6)=1.32936292837426009323670176765+50.1158017168222801180986619468 \text{ } \mathrm{I} \\
&\mathrm{x}(7)=1.32936292837426009323670176765+58.4684353362926601377817722712 \text{ } \mathrm{I} \\
&\mathrm{x}(8)=1.32936292837426009323670176765+66.8210689557630401574648825957 \text{ } \mathrm{I} \\
&\mathrm{x}(9)=1.32936292837426009323670176765+75.1737025752334201771479929202 \text{ } \mathrm{I} \\
&\mathrm{x}(10)=1.32936292837426009323670176765+83.5263361947038001968311032446 \text{ } \mathrm{I} \\
&\mathrm{x}(11)=1.32936292837426009323670176765+91.8789698141741802165142135691 \text{ } \mathrm{I} \\
&\mathrm{x}(12)=1.32936292837426009323670176765+100.231603433644560236197323894 \text{ } \mathrm{I} \\
&\mathrm{x}(13)=1.32936292837426009323670176765+108.584237053114940255880434218 \text{ } \mathrm{I}
\end{aligned}$\\\\\\
Infinity Roots\\\\\\
\textbf{6.2.Categories of complex roots corresponding to the case b:}\\\\
The basics Programs in mathematica.\\\\
\textbf{Constants:} ${m_{1}}=\frac{a}{c}=\sqrt[5]{43} ; m_{2}=\frac{b}{c}=-\sqrt[3]{31}$;\\[2pt]
$\mathrm{d}:=15 ; \mathrm{p}:=\pi ; \mathrm{q}:=\mathrm{e};$\\[2pt]
$\mathrm{h} 1\left[\mathrm{x}_{-}\right]:=(\mathrm{m} 1)^{\wedge} \mathrm{x}+\mathrm{p}^{*}(\mathrm{m} 2)^{\wedge} \mathrm{x}-\mathrm{q}$;\\[2pt]
$\mathrm{h 2\left[\mathrm{x}_{-}\right]:=-1/p^{*} (\mathrm{m} 1)}^{\wedge} \mathrm{x}+\mathrm{q}/\mathrm{p};$\\[2pt]
For $[\mathrm{k}=-\mathrm{d}, \mathrm{k}<=\mathrm{d}, \mathrm{k}++$,\\[2pt]
$\mathrm{f}\left[\mathrm{u}_{-}\right]:=\left(\log [\mathrm{u}]+\left(2^{*} \mathrm{k}^{*} \pi^{*} \mathrm{I}\right)\right) / \log [\mathrm{m} 2]$;\\[2pt]
$x r=$ Nest $[\mathrm{f}[(\mathrm{h} 2[\#])] \&, 0,10]$;\\[2pt]
$\mathrm{FQ} 2=\mathrm{N}[\mathrm{xr}, 10]$;\\[2pt]
Print["Approach $=$ ", $\mathrm{N}[\mathrm{h} 1[\mathrm{y}], 30]/$. FindRoot $[\mathrm{h} 1[\mathrm{y}]==0,\{\mathrm{y}, \mathrm{FQ} 2\}$, WorkingPrecision$\to 30], " \mathrm{x}(", \mathrm{k}, ")=", \mathrm{~N}[\mathrm{z}] /$.FindRoot $[\mathrm{h} 1[\mathrm{z}]==0,\{\mathrm{z}, \mathrm{FQ} 2\}$,
WorkingPrecision $\to 30]]]$
\newpage\noindent
\textbf{For the first 21 roots (from infinite complex Roots) we get.}\\\\
\textbf{Results 9:}\\\\
$\begin{aligned}
&\mathrm{x}(-10)=-17.6708670992525382500210144891-6.39244149043165467954365963797 \text{ } \mathrm{I}\\
&\mathrm{x}(-9)=-15.9052625797348216939204909555-5.74913079189852942523203026835 \text{ } \mathrm{I}\\
&\mathrm{x}(-8)=-14.1396581085817360423637483024-5.10582168159023406815737164189 \text{ } \mathrm{I}\\
&\mathrm{x}(-7)=-12.3740510102774718677920189344-4.46251796194974345001420000146 \text{ } \mathrm{I}\\
&\mathrm{x}(-6)=-10.6084256713969881208311412561-3.81922763811319342747954640732 \text{ } \mathrm{I}\\
&\mathrm{x}(-5)=-8.84271588555983753545329276412-3.17595005668124265151716787800 \text{ } \mathrm{I}\\
&\mathrm{x}(-4)=-7.07670147668304967743649663054-2.53256705974439390398062089406 \text{ } \mathrm{I}\\
&\mathrm{x}(-3)=-5.30984992287635303352948277895-1.88829785677614216529548120143 \text{ } \mathrm{I}\\
&\mathrm{x}(-2)=-3.54170969394503471525566086214-1.23955079654876832102220047884 \text{ } \mathrm{I}\\
&\mathrm{x}(-1)=-1.77731279444547086200298016750-0.572568978830005695936055265763 \text{ } \mathrm{I}\\
&\mathrm{x}(0)=-0.0751212877365536830355891729551+0.152969515224423777278081181818 \text{ } \mathrm{I}\\
&\mathrm{x}(1)=1.32451209943781285438925570727+0.725825584824742783606616383329 \text{ } \mathrm{I}\\
&\mathrm{x}(2)=2.96729579227675451503974652116+0.814278965513704256780960100130 \text{ } \mathrm{I}\\
&\mathrm{x}(3)=4.89173812645201749253290815756+0.991485702541660584621184657797 \text{ } \mathrm{I}\\
&\mathrm{x}(4)=6.85122399379713247579236902591+1.22319558394859143399479007723 \text{ } \mathrm{I}\\
&\mathrm{x}(5)=8.81785390561844772357282670808+1.46634612658618689421058998067 \text{ } \mathrm{I}\\
&\mathrm{x}(6)=10.7864240033118422313395583168+1.71180405846075920210445681490 \text{ } \mathrm{I}\\
&\mathrm{x}(7)=12.7555184144055098977865402109+1.95770093244574709918109959409 \text{ } \mathrm{I}\\
&\mathrm{x}(8)=14.7247481449571241878636245713+2.20367427809457833024642873265 \text{ } \mathrm{I}\\
&\mathrm{x}(9)=16.6940113070916808934037897543+2.44965907041596002761446021192 \text{ } \mathrm{I}\\
&\mathrm{x}(10)=18.6632824183316702903725095336+2.69564502316182224459278506286 \text{ } \mathrm{I}
\end{aligned}$\\\\\\
Infinity Roots\\\\\\
\textbf{IV.7. Solving the equation $\boldsymbol{\zeta(x)=x, x \in C}$, in the set of complexes. [1,2,6,7,8]}\\\\
To solve the equation we use of the 2 functional equations of Riemann. We use 4 special algorithms per case. Riemann's 2 functional equations are a function of $\zeta(), \cos (), \sin (),(2 \pi)^{x}$, Gamma( $)$ and depending on each function they produce the fields of the roots.
\[\begin{aligned}
&\frac{\zeta(x)}{\zeta(1-x)}-2(2 \pi)^{x-1} \cdot \sin \left(\frac{\pi x}{2}\right) \cdot \operatorname{Gamma}(1-x) \\
&\frac{\zeta(1-x)}{\zeta(x)}-2(2 \pi)^{-x} \cdot \cos \left(\frac{\pi x}{2}\right) \cdot \operatorname{Gamma}(x)
\end{aligned}\]
The functions that represent this method are cases resulting in detail per case:\\\\
We define the functions (cases a \& b)\\
\[f(x)=\frac{\zeta(x)}{\zeta(1-x)}-2(2 \pi)^{x-1} \cdot\sin \left(\frac{\pi x}{2}\right) \cdot \operatorname{Gamma}(1-x), \text{ } \zeta(x)=x\]\\
Case a: for $k=1, f_{1}(x)=\mathrm{x} / \zeta(1-\mathrm{x}) \cdot 1 / 2 \cdot(2 \cdot \pi)^{\wedge}(1-\mathrm{x}) / \operatorname{Gamma}(1-\mathrm{x})$.\\\\\\
$\left.a_{1}\right)$\\[4pt]
i) $\therefore \Phi_{1, \mu}\left(u_{1}\right)=\left(2\left(\pi-\arcsin \left(\mathrm{u}_{1}\right)+2 \cdot \pi \cdot \mathrm{w}\right)\right) / \pi, w \in Z.$\\[2pt]
ii) $\therefore{ }_{\sigma_{1}}^{m_{2}} x_{\#+1}^{\mu}=\Phi_{k, \mu}\left(\frac{1}{a_{k}}\left(-\sum_{i>1}^{2} a_{i} \sigma_{i}(\#)-a_{0}\right)\right)=\Phi_{1, \mu}\left(f_{1}(\#)\right)$.\\\\
$\left.a_{2}\right)$\\[4pt]
i) $\therefore \Phi_{2, \mu}\left(u_{2}\right)=\left(2\left(\arcsin \left(\mathrm{u}_{2}\right)+2 \cdot \pi \cdot \mathrm{w}\right)\right) / \pi, w \in Z.$\\[2pt]
ii) $\therefore{ }_{\sigma_{1}}^{m_{2}} x_{\#+1}^{\mu}=\Phi_{k, \mu}\left(\frac{1}{a_{k}}\left(-\sum_{i>1}^{2} a_{i} \sigma_{i}(\#)-a_{0}\right)\right)=\Phi_{1, \mu}\left(f_{1}(\#)\right)$.\\\\\\
I am using the function $f(x)=\frac{\operatorname{Gamma}\left(\frac{x}{2}\right)}{\operatorname{Gamma}\left(\frac{1-x}{2}\right)}-\pi^{x-1 / 2} \cdot \frac{\zeta(1-x)}{\zeta(x)}$.\\\\
Case b: for $k=1$, $f_{2}(x)=\frac{\operatorname{Gamma}\left(\frac{x}{2}\right)}{\operatorname{Gamma}\left(\frac{1-x}{2}\right)} \cdot \frac{\zeta(x)}{\zeta(1-x)}$\\\\\\
We define the functions (cases a \& b)\\\\
$\left.b_{1}\right)$\\
i) $\therefore \Phi_{3, \mu}\left(u_{1}^{\prime}\right)=-1 / 2+\left(\log \left(\mathrm{u}_{1}^{\prime}\right)+2 \cdot \mathrm{w} \cdot \pi \cdot \mathrm{I}\right) / \log (\pi), w \in Z.$\\
ii) $\underset{\sigma_{2}^{\prime}}{m_{1}} x_{\#+1}^{\mu}=\Phi_{1, \mu}^{\prime}\left(f_{2}(\#)\right)$.\\
Case b: for $k=1, f_{2}(x)=\frac{\operatorname{Gamma}\left(\frac{x}{2}\right)}{\operatorname{Gamma}\left(\frac{1-x}{2}\right)} \cdot \frac{\zeta(x)}{\zeta(1-x)}$:\\
$\left.b_{2}\right)$\\
i) $\therefore \Phi_{4, \mu}\left(u_{2}^{\prime}\right)=1 / 2+\left(\log \left(\mathrm{u}_{2}^{\prime}\right)+2 \cdot \mathrm{w} \cdot \pi \cdot \mathrm{I}\right) / \log (\pi), w \in Z.$\\
ii) $\therefore{ }_{\sigma_{1}}^{m_{2}} x_{\#+1}^{\mu}=\Phi_{2, \mu}^{\prime}\left(f_{2}(\#)\right)$.\\\\\\
\textbf{Numerical calculations for Complex roots:}\\\\
\textbf{7.1. Category of complex roots corresponding to the case $\boldsymbol{\alpha(\alpha 1)}$:}\\\\
The basics Programs in mathematica.\\\\
Unprotect[G1];\\[2pt]
$\mathrm{f1}[x_{-}]:=\mathrm{x} / \operatorname{Zeta}[1-\mathrm{x}]^{*}1/2^{*}(2^{*} \pi)^{\wedge}(1-\mathrm{x}) / \operatorname{Gamma}[1-\mathrm{x}]$;\\[2pt]
$\mathrm{f}[x_{-}]:=\mathrm{x}-2^{*}\left(2^{*} \pi\right)^{\wedge}(\mathrm{x}-1)^{*} \operatorname{Sin}\left[\pi^{*} \mathrm{x} / 2\right]^{*} \operatorname{Gamma}[1-\mathrm{x}]^{*} \operatorname{Zeta}[1-\mathrm{x}]$;\\[2pt]
$\mathrm{h}[\mathrm{x_{-}}]:=\operatorname{Zeta}[\mathrm{x}]-\mathrm{x}$;\\[2pt]
Off[FindRoot::Istol,FindRoot::cvmit,General::stop,FindRoot::bddir,FindRoot::nInum,ReplaceAll::reps,\\[2pt]
Infinity::indet,FindRoot::srect,General::munfl];\\[2pt]
$\mathrm{G1}[d_{-},q_{-},t_{-},s_{-},s1_{-},m_{-}]:=$\\[2pt]
For $[\mathrm{k}=-\mathrm{d}, \mathrm{k}<=\mathrm{5}, \mathrm{k}++$,\\[2pt]
$\mathrm{F1}[u_{-}]:=N[(2 (\pi-\operatorname{ArcSin}[u]+2 \pi \mathrm{k})) / \pi]$;\\[2pt]
$\mathrm{x r}=\operatorname{Nest}[\mathrm{F1}[\mathrm{f 1}[\#]] \&, \mathrm{q}, \mathrm{t}]$;\\[2pt]
$\mathrm{FQ3}= \mathrm{N} [\mathrm{x r}, \mathrm{8 0}]$;\\[2pt]
$\mathrm{s3}=\mathrm{N}[\mathrm{y}] /$.FindRoot $[\mathrm{f}[\mathrm{y}]==\mathrm{0},\{\mathrm{y}, \mathrm{F Q 3}\}$, WorkingPrecision $\to\mathrm{m}]$\\[2pt]
$\mathrm{s} 4=\mathrm{N}[\mathrm{y} 1] / \text {.FindRoot }[\mathrm{f}[\mathrm{y} 1]==0,\{\mathrm{y}, \mathrm{s} 3\}, \text { WorkingPrecision}\to\mathrm{m}]$;\\[2pt]
$\text {If [Abs}[s 3-s 4]]<s \text { \&\& Abs[s4]}<s1,$\\[2pt]
$\text {Print ["Approach=", }[h[y 1], m] / . \text { FindRoot }[h[y 1]==0,\{y 1, s 4\}, \text { WorkingPrecision}\to\mathrm{m}],$\\[2pt]
$\left.\left.\left."x\left({ }^{\prime \prime}, \mathrm{k}, "\right)=", N[y 1] / .  \text{FindRoot}[h[y 1]==0,\{\mathrm{y}, s 4\}, \text {WorkingPrecision }\to\mathrm{m}]\right], \text{Loopback}\right]\right]$\\[2pt]
Protect[G1];\\[2pt]
$\mathrm{G} 1[20,1 / 2,10,0.00001,75,70]$\\\\
\textbf{Results 10:}\\\\
$\begin{aligned}
&x(-9)=-33.999999999683437233174912281325947 \\
&x(-8)=-29.999999800525511561320078724806606 \\
&x(-7)=-25.999927036563283322479328779219027 \\
&x(-6)=-21.985531995733483012809206667804255 \\
&x(-5)=-17.729166140171729124345484517377206-1.1113876772379532029008867533810019 \text{ } \mathrm{I}\\
&x(-4)=-14.613492351835987588071919186650501+3.1007918168141534354604138798653371 \text{ } \mathrm{I}\\
&x(-3)=-10.821148006023595165258300978174459-5.3193414062340658907809379747379031 \text{ } \mathrm{I}\\
&x(-2)=-5.2793950307843170303114858336796244+8.8027675740832589731484730202165154 \text{ } \mathrm{I}\\
&\mathrm{x}(-1)=-0.29590500557521395564723783108304803 \\
&\mathrm{x}(0)=-0.29590500557521395564723783108304803 \\
&x(1)=1.8337726516802713962456485894415236
\end{aligned}$\\\\\\
Infinity Roots\\\\\\
The \textbf{approach from top to bottom} is of the order of $10^{\wedge}(-22)$ to $10^{\wedge}(-36)$.\\\\\\
\textbf{7.2. Category of complex roots corresponding to the case a2:}\\\\
The relevant basics Programs \textbf{in mathematica.}\\\\
Unprotect[G1];\\[2pt]
f1$\left[\mathrm{x}_{-}\right]:=\mathrm{x}/$Zeta$[1-\mathrm{x}]^{*} 1 / 2^{*}(2^{*} \pi)^{\wedge}(1-\mathrm{x}) /$Gamma$[1-\mathrm{x}]$;\\[2pt]
$\mathrm{f}[\mathrm{x_{-}}]:=\mathrm{x}-2^{*}(2^{*} \pi)^{\wedge}(\mathrm{x}-1)^{*} \operatorname{Sin}\left[\pi^{*} \mathrm{x} / 2\right]^{*}\operatorname{Gamma}[1-\mathrm{x}]^{*} \operatorname{Zeta}[1-\mathrm{x}]$;\\[2pt]
$\mathrm{h}\left[\mathrm{x}_{-}\right]:=$Zeta$[\mathrm{x}]-\mathrm{x}$;\\[2pt]
$G1[d_{-},q_{-},t_{-},s_{-},s1_{-},m_{-}]:=$\\[2pt]
For$[\mathrm{k}=-\mathrm{d}, \mathrm{k}<=\mathrm{2}, \mathrm{k}++$,\\[2pt]
F2$\left.\left[u_{-}\right]:=N[(2(\operatorname{ArcSin}[u]+2 \pi k)) / \pi\right]$;\\[2pt]
$\mathrm{x r}=\operatorname{Nest}[\mathrm{F} 2[\mathrm{f} 1[\#]] \mathcal{\&}, \mathrm{q}, \mathrm{t}]$;\\[2pt]
FQ3 =N[xr,80];\\[2pt]
s3=N[y]/.FindRoot$[\mathrm{f}[\mathrm{y}]==0,\{\mathrm{y}, \mathrm{F Q} 3\}$, WorkingPrecision $\to\mathrm{m}]$;\\[2pt]
$\mathrm{s} 4=\mathrm{N}[\mathrm{y} 1] /$.FindRoot $[\mathrm{f}[\mathrm{y} 1]==\mathrm{0},\{\mathrm{y} 1, \mathrm{s} 3\}$, WorkingPrecision$\to$m];\\[2pt]
If $[\mathrm{Abs}[\mathrm{s} 4]<\mathrm{s 1}$,Print ["x(", k,")=",N[y1]/.FindRoot $[h[y 1]==0,\{\mathrm{y 1}, \mathrm{s} 4\}$,\\[2pt]
WorkingPrecision$\to$m]],Loopback]]\\[2pt]
Protect[G1];\\[2pt]
G1$[20,1 / 2,10,0.00001,75,70]$
\newpage\noindent
\textbf{Results 11:}\\\\
$\begin{aligned}
&x(-8)=-32.000000008467671636378705967757058 \\
&x(-7)=-28.000004102805105685691855067848468 \\
&x(-6)=-24.001107180111899006547352622573465 \\
&x(-5)=-20.131186308008286295709692618189034 \\
&x(-4)=-14.613492351835987588071919186650501+3.1007918168141534354604138798653371 \text{ } \mathrm{I}\\
&x(-3)=-10.821148006023595165258300978174459+5.3193414062340658907809379747379031 \text{ } \mathrm{I}\\
&x(-2)=1.8337726516802713962456485894415236 \\
&x(-1)=-5.2793950307843170303114858336796244+8.8027675740832589731484730202165154 \text{ } \mathrm{I}\\
&x(1)=-0.29590500557521395564723783108304803 \\
&x(2)=-0.29590500557521395564723783108304803
\end{aligned}$\\\\\\
The \textbf{approach from top to bottom} is of the order of $10^{\wedge}(-22)$ to $10^{\wedge}(-36)$.In short, we have 4 complex in this area(categories)and the rest are real.\\\\\\
\textbf{7.3. Category of complex roots corresponding to the case b1:}\\\\
\textbf{The relevant basics Programs in mathematica}\\\\
Unprotect $[\mathrm{G} 1]$;\\[2pt]
ClearAll;\\[2pt]
$\mathrm{f2}[\mathrm{x_{-}}]:=\mathrm{x} / \operatorname{Gamma}[(1-\mathrm{x}) / 2]^{*} \mathrm{Gamma}[\mathrm{x} / 2] /(\mathrm{Zeta}[1-\mathrm{x}])$;\\[2pt]
$\mathrm{f}[\mathrm{x_{-}}]:=\mathrm{Gamma}[\mathrm{x} / 2]-(\mathrm{Zeta}[1-\mathrm{x}] / \mathrm{x})^{*} \pi^{\wedge}(\mathrm{x}-1 / 2)^{*} \operatorname{Gamma}[(1-\mathrm{x}) / 2]$;\\[2pt]
$\mathrm{h}\left[\mathrm{x}_{-}\right]:=$Zeta $[\mathrm{x}]-\mathrm{x}$;\\[2pt]
Off$[\text{FindRoot:: }$stol,FindRoot::cvmit,General::stop,FindRoot::bddir,FindRoot::nlnum,\\[2pt]
$\text{ReplaceAll::reps,Infinity::indet,FindRoot::srect,Less::nord,Divide::infy,Set::setraw];}$\\[2pt]
$G1[d_{-},q_{-},t_{-},s_{-},s1_{-},m_{-}]:=$\\[2pt]
For $[\mathrm{k}=-\mathrm{d}, \mathrm{k}<=\mathrm{d}, \mathrm{k}++$,\\[2pt]
$\mathrm{F} 3\left[\mathrm{u}_{-}\right]:=1 / 2+\left(\log [\mathrm{u}]+2^{*} \mathrm{k}^{*} \pi^{*} \mathrm{I}\right) / \log [\pi]$;\\[2pt]
$\mathrm{xr}=\mathrm{Nest}[\mathrm{F} 3[\mathrm{f} 2[\#]] \&, \mathrm{q}, \mathrm{t}]$;\\[2pt]
$\mathrm{FQ} 2=\mathrm{N}[\mathrm{xr}, 10] ; \mathrm{s} 3=0 ; \mathrm{s} 4=0$;\\[2pt]
$\mathrm{s} 3=\mathrm{N}[\mathrm{y} 4] /$.FindRoot $[\mathrm{f}[\mathrm{y} 4]==0,\{\mathrm{y} 4, \mathrm{FQ} 2\}$, WorkingPrecision$\to30]$;\\[2pt]
$\mathrm{s} 4=\mathrm{N}[\mathrm{y} 1] /$.FindRoot $[\mathrm{f}[\mathrm{y} 1]==0,\{\mathrm{y} 1, \mathrm{~s} 3\}$, WorkingPrecision$\to30]$;\\[2pt]
If $[\mathrm{Abs}[\mathrm{s} 4]<\mathrm{s} 1, \operatorname{Print}[" A p p r o a c h=", \mathrm{~N}[\mathrm{~h}[\mathrm{y}], \mathrm{m}] /$ FindRoot $[\mathrm{h}[\mathrm{y}]==0,\{\mathrm{y}, \mathrm{s} 4\}$,\\[2pt]
WorkingPrecision$\to\mathrm{m}], " \mathrm{x}(", \mathrm{k}, ")=", \mathrm{N}[\mathrm{y}] /$.FindRoot $[\mathrm{h}[\mathrm{y}]==0,\{\mathrm{y}, \mathrm{s} 4\}$,\\[2pt]
WorkingPrecision$\to\mathrm{m}], \mathrm{m}]] ;$ Loopback $]] / /$ Timing\\[2pt]
Protect[G1];\\\\
$\mathrm{G} 1[20,1 / 2,10,0.00001,110,30]$
\newpage\noindent
\textbf{Results 12:}\\\\
$\begin{aligned}
&x(-14)=-1.26139081555873261406489769808-87.7583685099246485627976563544 \text{ } \mathrm{I}\\
&\mathrm{x}(-13)=-1.12616495962146742040656000337-80.4582743596358549469056379024 \text{ } \mathrm{I}\\
&\mathrm{x}(-12)=-1.29290215735362288709093453220-75.6417409652417064465387628638 \text{ } \mathrm{I}\\
&x(-11)=-1.26227506725207853187875606381-70.3556950691515434161289031984 \text{ } \mathrm{I}\\
&x(-10)=-1.29443568561644493920364442108-65.2301042074151417066393823324 \text{ } \mathrm{I}\\
&\mathrm{x}(-9)=-1.41821905056698762133825325680-59.6454164781995695989688716976 \text{ } \mathrm{I}\\
&x(-8)=-1.28703179907482141334672081309-53.9110489785536948914586838211 \text{ } \mathrm{I}\\
&\mathrm{x}(-7)=-1.46547485014965062526389943588-48.0489077335066889394342997498 \text{ } \mathrm{I}\\
&x(-6)=-1.55341651249244969842510775727-41.5258170026898503231006830279 \text{ } \mathrm{I}\\
&x(-5)=-1.51767569589846950803032259371-38.1936599169588114867812150952 \text{ } \mathrm{I}\\
&x(-4)=-1.74955435401611051766624316632-30.8158449441671134562051307077 \text{ } \mathrm{I}\\
&x(-3)=-1.6934543711272871415911807768-26.5282868639131487782281764409 \text{ } \mathrm{I}\\
&x(-2)=-2.03689852362822026102533100727-21.9930698284317736851003882855 \text{ } \mathrm{I}\\
&x(-1)=-1.74955435401611051766624316629-30.8158449441671134562051307077 \text{ } \mathrm{I}\\
&x(0)=1.83377265168027139624564858944\\
&x(1)=-1.74955435401611051766624316629+30.8158449441671134562051307077 \text{ } \mathrm{I}\\
&x(2)=-2.03689852362822026102533100727+21.9930698284317736851003882855 \text{ } \mathrm{I}\\
&\mathrm{x}(3)=-1.6934543711272871415911807768+26.5282868639131487782281764409 \text{ } \mathrm{I}\\
&\mathrm{x}(4)=-1.74955435401611051766624316632+30.8158449441671134562051307077 \text{ } \mathrm{I}\\
&\mathrm{x}(5)=-1.51767569589846950803032259371+38.1936599169588114867812150952 \text{ } \mathrm{I}\\
&x(6)=-1.55341651249244969842510775727+41.5258170026898503231006830279 \text{ } \mathrm{I}\\
&x(7)=-1.46547485014965062526389943588+48.0489077335066889394342997498 \text{ } \mathrm{I}\\
&\mathrm{x}(8)=-1.28703179907482141334672081309+53.9110489785536948914586838211 \text{ } \mathrm{I}\\
&x(9)=-1.41821905056698762133825325680+59.6454164781995695989688716976 \text{ } \mathrm{I}\\
&x(10)=-1.29443568561644493920364442108+65.2301042074151417066393823324 \text{ } \mathrm{I}\\
&x(11)=-1.26227506725207853187875606381+70.3556950691515434161289031984 \text{ } \mathrm{I}\\
&x(12)=-1.29290215735362288709093453220+75.6417409652417064465387628638 \text{ } \mathrm{I}\\
&x(13)=-1.12616495962146742040656000337+80.4582743596358549469056379024 \text{ } \mathrm{I}\\
&\mathrm{x}(14)=-1.26139081555873261406489769808+87.7583685099246485627976563544 \text{ } \mathrm{I}\\
\end{aligned}$\\\\
Infinity Roots\\\\
\textbf{7.4. Category of complex roots corresponding to the case b2:}\\\\
\textbf{The relevant basics Programs in mathematica}\\\\
Unprotect$[\mathrm{G} 1]$; ClearAll;\\[2pt]
$\mathrm{f} 2\left[\mathrm{x}_{-}\right]:=\mathrm{x} / \operatorname{Gamma}[(1-\mathrm{x}) / 2]^{*} \operatorname{Gamma}[\mathrm{x}/2]/($Zeta$[1-\mathrm{x}])$;\\[2pt]
$\mathrm{f}[\mathrm{x_{-}}]:=\operatorname{Gamma}[\mathrm{x} / 2]-(\operatorname{Zeta}[1-\mathrm{x}] / \mathrm{x})^{*} \pi^{\wedge}(\mathrm{x}-1 / 2)^{*} \operatorname{Gamma}[(1-\mathrm{x})/2]$;\\[2pt]
$\mathrm{h}\left[\mathrm{x}_{-}\right]:=\mathrm{Zeta}[\mathrm{x}]-\mathrm{x}$;\\[2pt]
Off [FindRoot::lstol,FindRoot::cvmit,General::stop,FindRoot::bddir,FindRoot::nlnum,\\[2pt]
ReplaceAll::reps,Infinity::indet,FindRoot::srect,Less::nord,Divide::infy,Set::setraw];\\[2pt]
$\mathrm{G} 1\left[\mathrm{d}_{-}, \mathrm{q}_{-}, \mathrm{t}_{-}, \mathrm{s}_{-}, \mathrm{s} 1_{-}, \mathrm{m}_{-}\right]:=$\\[2pt]
For $[\mathrm{k}=-\mathrm{d}, \mathrm{k}<=\mathrm{d}, \mathrm{k}++$,\\[2pt]
$\mathrm{F} 4\left[\mathrm{u}_{-}\right]:=-1 / 2+\left(\log [\mathrm{u}]+2^{*} \mathrm{k}^{*} \pi^{*} \mathrm{I}\right) / \log [\pi]$;\\[2pt]
$\mathrm{xr}=\mathrm{Nest}[\mathrm{F} 4[\mathrm{f} 2[\#]] \&, \mathrm{q}, \mathrm{t}] ; \mathrm{FQ} 2=\mathrm{N}[\mathrm{xr}, 10] ; \mathrm{s} 3=0 ; \mathrm{s} 4=0$;\\[2pt]
$\mathrm{s} 3=\mathrm{N}[\mathrm{y} 4]/$.FindRoot $[\mathrm{f}[\mathrm{y} 4]==0,\{\mathrm{y} 4, \mathrm{FQ} 2\}$, WorkingPrecision $\to30]$;\\[2pt]
$\mathrm{s} 4=\mathrm{N}[\mathrm{y} 1]/$.FindRoot $[\mathrm{f}[\mathrm{y} 1]==0,\{\mathrm{y} 1, \mathrm{~s} 3\}$, WorkingPrecision$\to30]$;\\[2pt]
If $[\operatorname{Abs}[\mathrm{s} 4]<\mathrm{s} 1, \operatorname{Print}[" \operatorname{Approach}=", \mathrm{~N}[\mathrm{~h}[\mathrm{y}], \mathrm{m}] / .$ FindRoot $[\mathrm{h}[\mathrm{y}]==0,\{\mathrm{y}, \mathrm{s} 4\}$,\\[2pt]
WorkingPrecision$\to\mathrm{m}], " \mathrm{x}(", \mathrm{k}, ")=", \mathrm{zl}=\mathrm{N}[\mathrm{y}] / .$ FindRoot $[\mathrm{h}[\mathrm{y}]==0,\{\mathrm{y}, \mathrm{s} 4\}$,\\[2pt]
WorkingPrecision$\to\mathrm{m}]]$; Loopback$]]//$Timing\\[2pt]
Protect[G1];\\[2pt]
$\mathrm{G} 1[20,1 / 2,10,0.00001,110,30]$\\\\
\textbf{Results 13:}\\\\
$\begin{aligned}
&x(-14)=-1.23852071108706795556935851810-85.4090090247636352615544742032 \text{ } \mathrm{I}\\
&x(-13)=-1.29918130011957305015490858296-77.9865709309247937726779884803 \text{ } \mathrm{I}\\
&x(-12)=-1.15953302973918715196482277792-72.9767004265557497382345340656 \text{ } \mathrm{I}\\
&x(-11)=-1.33095139800651049656449887081-67.7834609908689492473936882734 \text{ } \mathrm{I}\\
&x(-10)=-1.1924742451865064471089288392-62.3193286218473504568063554479 \text{ } \mathrm{I}\\
&x(-9)=-1.3366144881440349550167977808-56.9002016278718488305446496687 \text{ } \mathrm{I}\\
&x(-8)=-1.40616347663699362699926390436-50.9143847707163058133311938748 \text{ } \mathrm{I}\\
&x(-7)=-1.32728372617780452008783857748-44.7148157580246341928779393786 \text{ } \mathrm{I}\\
&x(-6)=-1.55341651249244969842510775728-41.5258170026898503231006830279 \text{ } \mathrm{I}\\
&x(-5)=-1.51383649144589060970515941062-34.4539606289612051139663453008 \text{ } \mathrm{I}\\
&x(-4)=-1.74955435401611051766624316641-30.8158449441671134562051307077 \text{ } \mathrm{I}\\
&x(-3)=-2.03689852362822026102533100727-21.9930698284317736851003882855 \text{ } \mathrm{I}\\
&\mathrm{x}(-2)=-2.38593568342424543039818327574-16.2709870621967901723315117424 \text{ } \mathrm{I}\\
&x(-1)=-1.5138364914458906097051594106-34.4539606289612051139663453008 \text{ } \mathrm{I}\\
&x(0)=-1.55341651249244969842510775736+41.5258170026898503231006830279 \text{ } \mathrm{I}\\
&x(1)=-1.5138364914458906097051594106+34.4539606289612051139663453008 \text{ } \mathrm{I}\\
\end{aligned}$\\
$\begin{aligned}
&x(2)=-2.38593568342424543039818327574+16.2709870621967901723315117424 \text{ } \mathrm{I}\\
&x(3)=-2.03689852362822026102533100727+21.9930698284317736851003882855 \text{ } \mathrm{I}\\
&x(4)=-1.74955435401611051766624316641+30.8158449441671134562051307077 \text{ } \mathrm{I}\\
&x(5)=-1.51383649144589060970515941062+34.4539606289612051139663453008 \text{ } \mathrm{I}\\
&x(6)=-1.55341651249244969842510775728+41.5258170026898503231006830279 \text{ } \mathrm{I}\\
&x(7)=-1.32728372617780452008783857748+44.7148157580246341928779393786 \text{ } \mathrm{I}\\
&x(8)=-1.40616347663699362699926390436+50.9143847707163058133311938748 \text{ } \mathrm{I}\\
&x(9)=-1.3366144881440349550167977808+56.9002016278718488305446496687 \text{ } \mathrm{I}\\
&x(10)=-1.1924742451865064471089288392+62.3193286218473504568063554479 \text{ } \mathrm{I}\\
&x(11)=-1.33095139800651049656449887081+67.7834609908689492473936882734 \text{ } \mathrm{I}\\
&x(12)=-1.15953302973918715196482277792+72.9767004265557497382345340656 \text{ } \mathrm{I}\\
&x(13)=-1.29918130011957305015490858296+77.9865709309247937726779884803 \text{ } \mathrm{I}\\
&x(14)=-1.23852071108706795556935851810+85.4090090247636352615544742032 \text{ } \mathrm{I}\\
\end{aligned}$\\\\
Infinity Roots\\\\
The approach from top to bottom is of the order of $10^{\wedge}(-25)$ to $10^{\wedge}(-30)$.
\newpage\noindent
\textbf{Epilogue}\\\\
As we can see, each solution of an equation must be organized by categorizing its roots, so that the full range of roots is covered and not individually at some intervals. It is not always easy to choose the initial value, but we know 2 regions the interval $(-1,1)$ and the maximum value of the full function which is into region $[1, r>>1$ ). The choice and arrangement of examples is due to colleague Etairi Biragova thanks.
\vspace{2\baselineskip}
\phantomsection
\thispagestyle{plain}
\renewcommand\refname{References}

\end{document}